\documentclass[reqno,11pt]{amsart}
\usepackage[letterpaper,hmargin=1.25in,vmargin=0.9in]{geometry}
\usepackage[T1]{fontenc}
\usepackage{lmodern}
\usepackage{microtype}
\usepackage{tikz}
\usepackage{mathtools,amssymb,mathrsfs}
\usepackage[unicode,hidelinks]{hyperref}
\allowdisplaybreaks[1]
\newif\ifmarkdefinitions
\markdefinitionstrue
\definecolor{termink}{HTML}{154E73}
\definecolor{termborder}{HTML}{3981A6}
\definecolor{termfill}{HTML}{E9F5FA}
\ifmarkdefinitions
  \newcommand{\reviewterm}[1]{%
    \begingroup
    \setlength{\fboxsep}{1.1pt}%
    \setlength{\fboxrule}{0.45pt}%
    \fcolorbox{termborder}{termfill}{\textcolor{termink}{\upshape #1}}%
    \endgroup}
\else
  \newcommand{\reviewterm}[1]{#1}
\fi

\newcommand{\R}{\mathbb{R}}
\newcommand{\Z}{\mathbb{Z}}

\DeclareMathOperator{\conv}{conv}

\DeclareMathOperator{\interior}{int}

\DeclareMathOperator{\vol}{vol}
\DeclareMathOperator{\covol}{covol}

\DeclareMathOperator{\GL}{GL}
\newcommand{\doi}[1]{DOI: \href{https://doi.org/#1}{\nolinkurl{#1}}}

\theoremstyle{plain}
\newtheorem{theorem}{Theorem}[section]
\newtheorem{proposition}[theorem]{Proposition}
\newtheorem{lemma}[theorem]{Lemma}
\newtheorem{corollary}[theorem]{Corollary}

\theoremstyle{definition}

\theoremstyle{remark}

\numberwithin{equation}{section}

\newcommand{\scarf}{\mathcal S}
\newcommand{\idx}{\operatorname{ind}}

\title[A covering criterion and tetrahedron coverings]{A criterion for lattice coverings and the lattice covering density of a tetrahedron}
\author{Fei Xue}
\address{School of Mathematical Sciences, Nanjing Normal University, Nanjing 210046, China}
\email{05429@njnu.edu.cn}
\date{\today}
\hypersetup{pdftitle={A criterion for lattice coverings and the lattice covering density of a tetrahedron},pdfauthor={Fei Xue},
pdfsubject={Discrete geometry},pdfkeywords={lattice covering, tetrahedron, covering density}}
\begin{document}
\begin{abstract}
We give a necessary and sufficient condition for the lattice translates
$K+L$ of a convex body to cover $\R^n$. The criterion uses finite
integer chains of oriented $n$-simplices with vertices in $L$, allowing
affinely dependent vertices. The translates cover precisely when such a
chain has zero boundary modulo $L$, has nonzero total signed volume,
and the vertices of each term lie in a translate of $-K$.
We apply this criterion to prove that the lattice covering density of a
tetrahedron is $125/63$, as conjectured by Dougherty and Faber.
For the standard tetrahedron, we use the lattice points on the facets of
maximal lattice-free tetrahedra to prove the sharp covolume bound $21/250$.
The lattices attaining this bound are exactly the images of the $84$-shape
lattice under the linear parts of the standard tetrahedron's affine symmetries.
\end{abstract}
\maketitle

\section{Introduction}
Let $K$ be a convex body in $\R^n$. Among the lattices $L$ for which
$K+L$ covers $\R^n$, how large can the covolume of $L$ be? Equivalently,
one asks for the \reviewterm{lattice covering density}
\[
       \theta^l(K)=\inf_{K+L=\R^n}\frac{\vol_n(K)}{\covol(L)},
\]
where $L$ ranges over lattices of full rank and $\covol(L)$ is the volume
of a fundamental parallelepiped of $L$. The ratio inside the infimum is
the density of the particular covering $K+L$. We give a criterion for a
fixed pair $(K,L)$ to be a covering and use it to determine the density
of a tetrahedron. By affine invariance it suffices to treat
$\Delta=\conv(0,e_1,e_2,e_3)$.

F\'ary~\cite{Fary} proved in 1950 that a triangle has lattice covering
density $3/2$, the largest value among planar convex bodies.
Januszewski~\cite{Jan} later proved that arbitrary translates of a triangle
give the same optimal density.

In a regular tetrahedron the edge midpoints span an octahedron of half
its volume. Truncating its six vertices gives a regular truncated
octahedron, a lattice tile in Fedorov's 1885 classification~\cite{Fed}.
The tile occupies $8/9$ of the octahedron and $4/9$ of the tetrahedron,
giving a tetrahedron covering of density $9/4$~\cite[p.~501]{DF}.
Aguil\'o, Fiol, and Garcia~\cite{AFG} used a discrete construction
based on the same lattice in their 1997 work on triple-loop networks
(see~\cite[p.~502]{DF}).
In 1999 Lassak~\cite[p.~456]{Lassak} conjectured that the tile has
maximum volume among all \reviewterm{zonotopes} (finite Minkowski sums of line
segments) in the tetrahedron.

A smaller density emerged from the degree--diameter problem for abelian
Cayley digraphs. In 1994 Hsu and Jia~\cite{HJ} studied cyclic groups
with three generators, one prescribed. Forcade and
Lamoreaux~\cite[p.~194]{FL} date the discovery of the $84$-shape to
the same year, crediting Dougherty and Faber and, independently,
Fiduccia, Zito, and Mann. Its $84$ unit cubes form a nonconvex lattice
tile inside $10\Delta$, giving density $\vol_3(10\Delta)/84=125/63<2$.
The tile occupies $63/125$ of the tetrahedron. By a theorem of
F\'ary and R\'edei~\cite[Satz~5, p.~208]{FR}, even the largest centrally
symmetric convex subset occupies only half.

Fiduccia, Forcade, and Zito~\cite{FFZ} gave an infinite family of
related digraphs, later recovered by dilation~\cite{AFP}. Forcade and
Lamoreaux~\cite{FL} proved local optimality of the $84$-shape in 2000.
In their 2004 paper, Dougherty and Faber~\cite[pp.~487, 501--507]{DF}
described how their search among cyclic groups improved on $9/4$ and
conjectured that $125/63$ is globally optimal. Cocke, Forcade, and Hall~\cite{CFH} subsequently studied
local optima within fixed combinatorial types of fundamental region.

For the global problem, Xue and Zong~\cite{XZ} obtained the lower bound
$1+2^{-16}$ in 2018. In 2023 Fu, Xue, and Zong~\cite{FXZ} used the
connection with abelian Cayley digraphs~\cite{DF,FFZ,KT} to raise the
lower bound to $25/18$, leaving the interval
\[
                 \frac{25}{18}\le \theta^l(\Delta)\le\frac{125}{63}.
\]
For arbitrary translative coverings of the tetrahedron, Li, Fu, and
Zhang~\cite{LFZ} proved the lower bound $1+\sqrt{2}/1152$.

\subsection{A lattice covering criterion}
Consider again a convex body $K\subset\R^n$, $n\ge1$, and a full-rank lattice
$L$ of covolume $D$. A finite nonempty set $F\subset L$ is
\reviewterm{admissible} if the corresponding translates have a common point. Equivalently,
\begin{equation}\label{eq:co-covering}
     \bigcap_{v\in F}(K+v)\ne\varnothing
       \quad\Longleftrightarrow\quad F\subset x-K\quad
       \text{for some }x\in\R^n .
\end{equation}
The common point $x$ may depend on $F$. An \reviewterm{admissible integer chain} is a
finite integer combination of oriented simplices on admissible vertex
sets, taken modulo lattice translations. Its \reviewterm{boundary} is the alternating
sum of the oriented faces of its terms. A \reviewterm{cycle} is a chain
$z$ with $\partial z=0$, so these faces cancel with multiplicity modulo $L$.

For $S=[v_0,\dots,v_n]$ put
\[
       d_L(S)=\frac{\det[v_1-v_0,\dots,v_n-v_0]}{D}
\]
and extend this definition linearly to chains. Even if the vertices of $S$ are affinely dependent, its faces
still enter the boundary, while $d_L(S)=0$. The following
criterion arose in joint work with Yanlu Lian.

\begin{theorem}[Covering criterion]\label{thm:closed-chain}
The translates $K+L$ cover $\R^n$ if and only if there is an admissible
integer $n$-chain $z$ such that $\partial z=0$ and $d_L(z)\ne0$.
When $K+L=\R^n$, one can choose $z$ with $d_L(z)=n!$.
\end{theorem}

Balitskiy~\cite[Lemma~5.3]{Bal} proved a sufficient condition for
covering. If $V$ is the vertex set of a triangulation of $\R^n$ and
each simplex can be translated into $-K$, every translate of $-K$
contains a point of $V$. For $V=L$ this gives $K+L=\R^n$.
For a lattice, Theorem~\ref{thm:closed-chain} gives the same conclusion
from any admissible cycle of nonzero degree, where the \reviewterm{degree} is
$d_L(z)/n!$. Its simplices may overlap, occur with either orientation,
or have zero volume. Conversely, from every lattice covering we construct
an admissible cycle of degree one. We use this orientation whenever we
choose a cycle from a covering. A related homological criterion for
coverage by small equal balls in the maximum metric on a flat torus is due to
Decreusefond et al.~\cite{DFRV}.

\subsection{The tetrahedron}
For the standard tetrahedron $\Delta$, the largest covolume is $21/250$.
Write $L_{84}$ for the corresponding \reviewterm{$84$-shape lattice},
\[
 L_{84}=\tfrac1{10}\{(x,y,z)\in\Z^3:35x+9y+2z\equiv0\pmod {84}\}.
\]
\begin{theorem}\label{thm:main}
For every tetrahedron $T$ of positive volume,
\[
                          \theta^l(T)=\frac{125}{63}.
\]
A lattice $L$ attains this density for $T$ if and only if
$T=t+A\Delta$ and $L=A L_{84}$ for some $t\in\R^3$ and
$A\in\GL_3(\R)$.
\end{theorem}
The dilation construction of Aguil\'o, Fiol, and P\'erez~\cite[Proposition~4.3]{AFP}
gives abelian Cayley digraphs on three generators with order $84k^3$
and diameter $10k-3$. The following bound shows that each has the
largest possible order among such digraphs of its diameter.

\begin{corollary}\label{cor:cayley}
A finite abelian group on three specified generators whose directed Cayley
graph has diameter $d$ has order at most $21(d+3)^3/250$. For every
positive integer $k$, equality is attained at $d=10k-3$.
\end{corollary}

In their work on integer programming, B\'ar\'any, Howe, and
Scarf~\cite{BHS} studied the periodic complex of maximal lattice-free
simplices with prescribed facet normals. For the normals of $-\Delta$,
assume that no nonzero vector of $L$ is parallel to a facet. In this
general position, each such tetrahedron has no
lattice point in its interior and one on each facet. These four points
label an abstract cell. Scarf's adjacency rule~\cite[Section~3]{Scarf}
describes how the contacts change across a common triangular face.
We apply Theorem~\ref{thm:closed-chain} to the degree-one cycle formed
by the cells modulo $L$ to show that the covering radius of $\Delta$
is the least $t$ for which every cell is admissible for $t\Delta$
(Corollary~\ref{cor:scarf-radius}).

Suppose a lattice gave a covering of smaller density than the
$84$-shape. By Lemma~\ref{lem:perturb} we could choose such a lattice
$L$ in general position. Must some cell have contacts $p_0,p_1,p_2,p_3$
for which $p_1-p_0,p_2-p_0,p_3-p_0$ form a basis of $L$?
Sections~\ref{sec:highindex}--\ref{sec:flat-components} give such a
cell. Section~\ref{sec:finite} uses this basis to bound the covolume
of $L$ and identify the lattices attaining the bound.

\subsection*{Acknowledgments}
The author thanks OpenAI's ChatGPT 6.2 Astra for helping to improve the
exposition and the proofs, and for designing the algorithms and writing
their descriptions. The author takes full responsibility for the
mathematical claims and conclusions.

In lattice coordinates, the volume bounds leave only finitely many possible
positions for neighbouring contact points. The programs enumerate candidate
contact configurations and verify the resulting determinant bounds in exact
arithmetic. The arXiv ancillary files contain the programs and data. To
reproduce all the checks, run \texttt{python3 verify.py} from the directory
containing the main ancillary \texttt{README.txt}. The individual commands
appear with the corresponding proofs.

\section{Lattice coverings and periodic chains}\label{sec:cycle-criterion}

\subsection{Two triangles in the plane}
Let $A_6(K)$ be the largest area of a centrally symmetric hexagon
contained in a planar convex body $K$, allowing a parallelogram as a
degenerate hexagon. F\'ary~\cite[Section~III, Lemma~1]{Fary} proved
that whenever $K+L$ covers the plane, $K$ contains such a polygon
which tiles by $L$. Its area is $\covol(L)$. Conversely, every
centrally symmetric hexagon or parallelogram tiles by a lattice.
Consequently
\[
                \theta^l(K)=\frac{\operatorname{area}(K)}{A_6(K)}.
\]

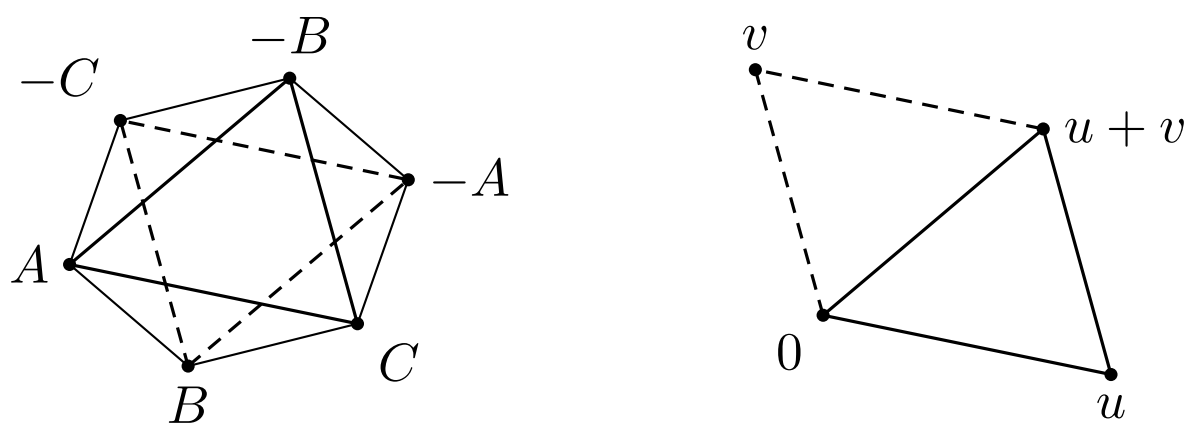
\begin{figure}[htb]
\centering
\begin{tikzpicture}[scale=1.12,line join=round,font=\small]
\begin{scope}
\coordinate (A) at (-1,-.25); \coordinate (B) at (-.3,-.85);
\coordinate (C) at (.7,-.6);  \coordinate (D) at (1,.25);
\coordinate (E) at (.3,.85);  \coordinate (F) at (-.7,.6);
\draw[thin] (A)--(B)--(C)--(D)--(E)--(F)--cycle;
\draw[semithick] (A)--(C)--(E)--cycle;
\draw[semithick,densely dashed] (B)--(D)--(F)--cycle;
\foreach \p/\position/\lab in
 {A/left/A,B/below/B,C/below right/C,D/right/-A,E/above/-B,F/above left/-C}
 {\fill (\p) circle (1.1pt) node[\position] {$\lab$};}
\end{scope}
\begin{scope}[shift={(3.45,-.55)}]
\coordinate (O) at (0,0); \coordinate (U) at (1.7,-.35);
\coordinate (V) at (-.4,1.45); \coordinate (W) at (1.3,1.1);
\draw[semithick,densely dashed] (O)--(V)--(W);
\draw[semithick] (O)--(U)--(W)--cycle;
\foreach \p/\position/\lab in
 {O/below left/0,U/below/u,V/above/v,W/right/{u+v}}
 {\fill (\p) circle (1.1pt) node[\position] {$\lab$};}
\end{scope}
\end{tikzpicture}
\caption{The alternating vertices of the hexagon give two triangles (left).
Translating them by $-A$ and $-B$ gives the two halves of a
fundamental parallelogram (right).}
\label{fig:planar-chain}
\end{figure}

Choose a nondegenerate centrally symmetric hexagon $H$ in $K$ and translate both so that
$H$ is centred at $0$. List its vertices
counterclockwise as $A,B,C,-A,-B,-C$. Put $u=C-A$ and $v=-B-C$.
The neighbour translations across the sides are, up to sign, $u,v,u+v$,
so $H$ tiles by the lattice $L_H=\Z u+\Z v$. The fundamental
parallelogram on $u,v$ is the sum, as an oriented chain, of
\[
              T_1=[0,u,u+v],\qquad T_2=[0,u+v,v].
\]
Within $H$, the alternating vertices give precisely the triangles
$A+T_1=[A,C,-B]$ and $B+T_2=[B,-A,-C]$. Their convex hull
is $H$ (Figure~\ref{fig:planar-chain}). Since $H=-H$, each of
$T_1,T_2$ can also be translated into $-K$. If $H$ is a parallelogram,
take two side vectors $u,v$ from one vertex in positive order. After
moving this vertex to $0$, the diagonal from $0$ to $u+v$ gives the
same two triangles for $L_H=\Z u+\Z v$.

On the lattice the common diagonal cancels in $\partial(T_1+T_2)$.
The outside edges cancel in pairs modulo $L_H$, and the signed
areas add to $\covol(L_H)$. In higher dimensions we retain
these two properties of the chain. Each simplex may still be
translated independently into $-K$.

\subsection{The degree of an admissible cycle}
Recall from~\eqref{eq:co-covering} that a finite nonempty set
$F\subset L$ is admissible if the translates $v+K$, $v\in F$, have
a common point. Regard each such $F=\{v_0,\ldots,v_k\}$ as an
abstract $k$-simplex, even when its lattice points are affinely
dependent. Its nonempty subsets are again admissible and give its
faces. These simplices form the \reviewterm{abstract simplicial complex} $N$.
Its \reviewterm{geometric realization} $|N|$ is obtained by taking a standard
simplex for each $F$ and gluing along common faces. A point in the
simplex on $F$ has \reviewterm{barycentric coordinates} $t_0,\ldots,t_k$, the
weights on its vertices satisfying $t_i\ge0$ and $\sum_i t_i=1$.

Translation by $w\in L$ takes the simplex on $F$ to that on $F+w$,
preserving barycentric coordinates. The \reviewterm{quotient} $X=|N|/L$ identifies
each point with all its lattice translates. All lattice vertices thus
give one vertex of $X$. Each orbit of simplices gives a cell of $X$,
with its faces read from any representative in $N$. There are only finitely
many translation classes of simplices, since moving one vertex to $0$
places all the others in the finite set $(K-K)\cap L$ by~\eqref{eq:co-covering}.
No nonzero translation preserves a finite nonempty vertex set. Fix a lattice basis
and list the vertices of each simplex in lexicographic order. This order
agrees on common faces and is preserved by translations. A $k$-chain on $X$
is then a sum $z=\sum_\sigma n_\sigma\sigma$ with one integer coefficient
for each translation class of $k$-simplices oriented in this way. In degree
$n$ these are precisely the admissible chains of Theorem~\ref{thm:closed-chain}.
The terms of the sum may be listed in any order. An odd permutation of the vertices of
one simplex reverses its orientation. If two faces of a simplex agree modulo $L$, both
occurrences enter its boundary.

The map $\widetilde a:|N|\to\R^n$ sends each vertex to its lattice
point and extends affinely over each simplex. Thus for a point $p$
with weights $t_i$ on the vertices $v_i$ we have
$\widetilde a(p)=\sum_{i=0}^k t_i v_i$.
Since $\widetilde a(p+w)=\widetilde a(p)+w$, this gives a map
$a:X\to T_L=\R^n/L$ by taking the image modulo $L$.

For an $n$-cycle $z$ in $X$, applying $a$ to each oriented simplex
of $z$ with its original coefficient gives a cycle on the torus. Its class is
$a_*[z]$, where $a_*$ is the map induced by $a$ on homology.
In Figure~\ref{fig:planar-chain}, for example, $T_1$ and $T_2$ form a
positively oriented fundamental parallelogram for $L_H$. Identifying
its opposite sides gives the torus $T_{L_H}$, so
$a_*[T_1+T_2]=[T_{L_H}]$, the fundamental class of the oriented
torus. In general, the degree of $z$ is the integer coefficient of
$[T_L]$ in its image class. In terms of signed volumes,
\begin{equation}\label{eq:cycle-degree}
              a_*[z]=\deg(z)[T_L],\qquad \deg(z)=d_L(z)/n!.
\end{equation}
Indeed, the volume form $dx_1\wedge\cdots\wedge dx_n/\covol(L)$
has integral one on $[T_L]$, which generates $H_n(T_L;\Z)$. Its
pullback integrates to $d_L(S)/n!$ on each oriented lattice simplex
$S$. Summing these integrals gives~\eqref{eq:cycle-degree}.
Thus degree one means that the signed volumes add to the volume of
one fundamental domain.

\subsection{From a cycle to a covering}
Let $z$ be an $n$-cycle in $X$ with $\deg(z)\ne0$. On the barycentric
subdivision of $N$ we construct a map $\widetilde g:|N|\to K+L$.
After passage to the quotient, its image of $z$ will represent the
nonzero class $a_*[z]$ in $T_L$. If the translates missed a point of
the torus, this class would have to come from the punctured torus,
whose top homology is zero.

\begin{proof}[Proof of sufficiency in Theorem~\ref{thm:closed-chain}]
For each admissible $F$ choose a point $x_F\in\bigcap_{v\in F}(v+K)$,
choosing once per translation orbit and setting $x_{F+w}=x_F+w$.
\ifmarkdefinitions\par\fi
The \reviewterm{barycentric subdivision} gives each nonempty face $F$ a vertex,
its barycentre $b_F$. Set $\widetilde g(b_F)=x_F$ and extend affinely
over each simplex of the subdivision. Such a simplex has vertices
$b_{F_0},\ldots,b_{F_k}$ for a chain $F_0\subsetneq\cdots\subsetneq F_k$.
Pick $v\in F_0$. Every $F_i$ contains $v$, so all the $x_{F_i}$ lie
in $v+K$. By convexity,
\[
       \widetilde g([b_{F_0},\ldots,b_{F_k}])
       =\conv(x_{F_0},\ldots,x_{F_k})\subset v+K.
\]
Thus $\widetilde g$ takes values in $K+L$ and commutes with
translations by $L$, giving a map $g:X\to T_L$ with image in $(K+L)/L$.

For $p\in|N|$, the straight segment from $\widetilde a(p)$ to
$\widetilde g(p)$ translates by $w$ when $p$ is replaced by $p+w$.
These segments therefore give a homotopy from $a$ to $g$ on taking
quotients. Homotopy invariance~\cite[Theorem~2.10]{Hatcher} gives
$g_*[z]=a_*[z]=\deg(z)[T_L]\ne0$.
If some $y\in T_L$ lay outside $(K+L)/L$, the map $g$ would factor
through $T_L\setminus\{y\}$. To see that $H_n(T_L\setminus\{y\};\Z)=0$,
centre a fundamental parallelepiped at a lift of $y$ and retract the
punctured parallelepiped radially onto its boundary. The retraction
fixes the boundary throughout and hence descends to the torus, where
the image of the boundary has dimension $n-1$. The factorization
would give $g_*[z]=0$, a contradiction. Hence $K+L=\R^n$.
\end{proof}

\subsection{From a covering to a cycle of degree one}
Suppose $K+L=\R^n$. We label the vertices of a sufficiently fine
$L$-periodic triangulation of $\R^n$ by nearest lattice points to
construct a map $\widetilde s:\R^n\to|N|$. The induced map
$s:T_L\to X$ sends the fundamental cycle of the triangulated torus to
the desired cycle $z$. Its degree will be one because $a\circ s$ is
homotopic to the identity on $T_L$.

Translating $K$ preserves the covering property and the admissible
vertex sets, so we may take $0\in\interior K$. Its \reviewterm{gauge}
$p_K(y)=\inf\{t\ge0:y\in tK\}$ measures the dilation needed to
reach $y$ in its own direction. For $x\in\R^n$ write
\[
 d_K(x,L)=\min_{v\in L}p_K(x-v),\qquad
 N_K(x)=\{v\in L:p_K(x-v)=d_K(x,L)\}.
\]
The points of $N_K(x)$ are the \reviewterm{nearest lattice points} to $x$ for this
gauge, which may take different values at $y$ and $-y$. On any bounded set of $x$, the bound $d_K(x,L)\le p_K(x)$
and the boundedness of $K$ leave only finitely many candidate lattice points. Thus $d_K(\cdot,L)$ is continuous and $L$-periodic. The \reviewterm{covering
radius} $\mu(K,L)$, the least $t$ for which $tK+L=\R^n$, is its
maximum over a closed fundamental parallelepiped.

\begin{proof}[Proof of necessity in Theorem~\ref{thm:closed-chain}]
To ensure that the labels on each simplex form an admissible set,
we first show that for some $\eta>0$, the lattice points chosen as
nearest to the vertices of any simplex of diameter less than $\eta$
are all nearest to a common point. Otherwise, translate a vertex of
each simplex in a sequence of counterexamples into a fixed closed
fundamental parallelepiped and let their diameters tend to zero.
Pass to a subsequence along which that vertex tends to $x$. All the
vertices then tend to $x$. Their chosen lattice points belong to a
fixed finite set, since the covering gives $d_K(p,L)\le1$ and hence
each nearest point to $p$ lies in $p-K$. Passing to a further
subsequence, the set of chosen labels is constant. For each such $v$,
continuity of $p_K$ and $d_K$ gives $p_K(x-v)=d_K(x,L)$. Thus every
chosen label belongs to $N_K(x)$, contrary to the choice of the simplices.

Take an $L$-periodic triangulation with each simplex of diameter less
than $\eta$, by subdividing the cubes of a sufficiently fine grid in
lattice coordinates according to coordinate order. At each vertex $p$
choose $\lambda(p)\in N_K(p)$, taking $\lambda(p+w)=\lambda(p)+w$ for
$w\in L$. The labels on each simplex are nearest to a common point
$x$. Since $d_K(x,L)\le1$, the translates $\lambda(p)+K$ for all
its vertices $p$ contain $x$. Thus the distinct labels span a simplex
of $N$. Sending $p$ to $\lambda(p)$ and extending over each simplex
gives the simplicial map $\widetilde s$. It commutes with lattice
translations and hence descends to $s:T_L\to X$.

Let $c$ be the sum of the positively oriented $n$-simplices of the
triangulated torus. Their faces cancel in pairs, so $c$ is a cycle
representing $[T_L]$. Its image chain $z$ is obtained by replacing
each vertex by its label, with a term counted as zero if a label is
repeated. The simplicial chain map commutes with boundary, so
$\partial z=0$ and $[z]=s_*[T_L]$. The composite $a\circ s$ lifts
to the piecewise affine map $\widetilde a\circ\widetilde s$ of $\R^n$
that sends each vertex $p$ to $\lambda(p)$. Its straight-line homotopy
to the identity commutes with translations and descends to $T_L$. Hence
\[
         a_*[z]=(a\circ s)_*[T_L]=[T_L].
\]
By~\eqref{eq:cycle-degree}, $\deg(z)=1$ and $d_L(z)=n!$.
\end{proof}

\section{Tetrahedron coverings and their lattice points}\label{sec:geometry}

We first exhibit degree-one cycles for the coverings of densities
$9/4$ and $125/63$. To rule out a covering of density below
$125/63$, we shall need to pass from a tetrahedron to its neighbour
across a prescribed triangle. Theorem~\ref{thm:closed-chain} requires
only that the triangular faces cancel with multiplicity. After passing
to general position, we shall choose a cycle in which each triangle
belongs to exactly two tetrahedra.

\subsection{Two covering lattices}
The covering of density $9/4$ admits a cycle of six tetrahedra.
Cutting off the four corners of a regular tetrahedron at its edge
midpoints leaves an octahedron of half its volume. A lattice covering
by this octahedron of density $9/8$ gives a tetrahedron covering of
density $9/4$~\cite[Section~8]{DF}. In coordinates for $6\Delta$,
a basis of the corresponding lattice $L$ consists of the three columns of
$4I-J$, where $J$ is the all-ones matrix. They are the permutations of
$(3,-1,-1)$. For each ordering $a,b,c$ of these columns take the
\reviewterm{Kuhn tetrahedron}
\[
                 [0,a,a+b,a+b+c].
\]
These six tetrahedra subdivide the fundamental parallelepiped. With
positive orientations their faces cancel modulo $L$ and their volumes
sum to $\covol(L)$, so they form a cycle of degree one. To check
admissibility, take $a=(3,-1,-1)$ and $b=(-1,3,-1)$. The coordinatewise
maximum of the four vertices is $h=(3,2,1)$. As $v$ runs through
$0,a,a+b,a+b+c$, the vectors $h-v$ are nonnegative and have coordinate
sums $6,5,4,3$, respectively. Thus $h$ lies in all four translates
$v+6\Delta$. Permuting coordinates gives a common point for each of the
other five terms. Theorem~\ref{thm:closed-chain} gives $6\Delta+L=\R^3$.
Since $\det(4I-J)=16$, the density is $\vol_3(6\Delta)/16=6^3/(6\cdot16)=9/4$.

\emph{A thirteen-term cycle at density $125/63$.}
For the $84$-shape lattice, take
\begin{equation}\label{eq:84-lattice}
       L=C\Z^3,\qquad C=\begin{pmatrix}2&1&-1\\2&-5&3\\-2&5&4\end{pmatrix}.
\end{equation}
Each column of $C$ satisfies the congruence $35x+9y+2z\equiv0\pmod{84}$
defining $10L_{84}$. Its solution lattice has index $84$, since
$35x+9y+2z$ takes the value $1$ at $(0,1,-4)$. As $|\det C|=84$,
we have $L=10L_{84}$.

The columns below are the fourteen lattice points $u_j\in C\Z^3$
used in the translates $u_j+10\Delta$. Write $[i,j,k]$ for
the oriented simplex $[0,u_i,u_j,u_k]$ with the indicated vertex order.
\[
\begin{array}{c|rrrrrrrrrrrrrr}
 j&1&2&3&4&5&6&7&8&9&10&11&12&13&14\\ \hline
 (u_j)_1&-2&-1&-4&-3&-6&-1&-5&-8&-3&-7&-2&-5&-1&-1\\
 (u_j)_2& 8& 5& 6& 3& 4& 3& 1& 2& 1&-1&-2&-1&-5&-7\\
 (u_j)_3&-1&-5& 1&-3& 3& 4&-1& 5& 6& 1& 2& 8&-2& 7
\end{array}
\]
Take the chain
\begin{equation}\label{eq:84-cycle}
\begin{aligned}
 z={}&[2,3,1]+[2,5,3]+[4,11,2]+[7,11,4]+[8,9,6]\\
     &+[10,11,7]+[11,5,2]+[11,6,5]+[11,8,6]\\
     &+[11,9,8]+[11,14,9]+[13,10,7]+[14,12,9].
\end{aligned}
\end{equation}
For the first term $[2,3,1]$ in~\eqref{eq:84-cycle}, the coordinatewise
maximum is $h=(0,8,1)$. The four vectors $h-v$ have coordinate sums
$9,10,6,4$ in vertex order, so $h$ belongs to every $v+10\Delta$.
Taking the coordinatewise maximum separately for each of the other
twelve terms gives nonnegative differences with coordinate sums at
most $10$ in every case. Thus every term is admissible for $10\Delta$.

To check closure, subtract
the least vertex in lexicographic order from each oriented triangular
face and sort its vertices, retaining the sign of the permutation.
The $52$ faces then give $26$ distinct triples, each with both signs.
For instance $[0,u_2,u_3]$ occurs with opposite signs in the first
two terms. The signed determinants of the thirteen terms, divided by $84$ and
listed as in~\eqref{eq:84-cycle}, are
\[
                1,1,0,0,1,0,1,1,-1,1,1,1,-1 .
\]
Thus $d_L(z)=8-2=6$, so $z$ has degree one. The three terms with
coplanar vertices have zero volume, but their faces enter the
cancellation. Applying Theorem~\ref{thm:closed-chain} to $z$ gives
$10\Delta+C\Z^3=\R^3$. Scaling by $1/10$ gives the covering
$\Delta+L_{84}$ of covolume $84/10^3=21/250$ and density $125/63$.
The ancillary command \texttt{python3 examples.py cycles} checks the faces,
signed determinants, and admissible scales for both examples.

For a lattice $B\Z^3$ with basis matrix $B$, the lower bound in
Theorem~\ref{thm:main} is equivalent to
\begin{equation}\label{eq:goal}
  \Delta+B\Z^3=\R^3\quad\Longrightarrow\quad |\det B|\le\frac{21}{250}.
\end{equation}
To prove~\eqref{eq:goal}, we first express the covering radius as a
minimum over degree-one cycles and use this description to reduce the
problem to lattices in general position.

\subsection{The radius of a finite set of lattice points}\label{sec:coordinates}
Fix $L=B\Z^3$. A lattice point with coordinates $p\in\Z^3$ in the
basis formed by the columns of $B$ has standard coordinates $Bp\in L$.
We call $p$ the \reviewterm{integer label} of $Bp$.
For a set $P\subset\Z^3$ of labels, write $BP=\{Bp:p\in P\}$
for the corresponding set of lattice points. Given a finite nonempty
$P$, what is the least $t\ge0$ for which the translates
$Bp+t\Delta$, $p\in P$, have a common point? By~\eqref{eq:co-covering}
this is the \reviewterm{radius}
\[
 r_B(P)=\min\{t\ge0:BP\subset x-t\Delta\text{ for some }x\in\R^3\}.
\]
A point $x$ lies in all the translates precisely when, for every $p\in P$,
the vector $x-Bp$ is coordinatewise nonnegative and its coordinate sum
is at most $t$. To express these four inequalities, write
\[
 g_i(p)=(Bp)_i\quad(1\le i\le3),\qquad
 g_4(p)=-\sum_{i=1}^3(Bp)_i,\qquad h_i(P)=\max_{p\in P}g_i(p).
\]
Each $h_i(P)$ is the position of the supporting plane of $BP$ in the
corresponding facet direction. The first three inequalities require
$x_i\ge h_i(P)$, and the fourth requires $\sum_{i=1}^3x_i+h_4(P)\le t$.
Among the points satisfying the first three, $h=(h_1(P),h_2(P),h_3(P))$
has the smallest coordinate sum. Thus $BP\subset h-r_B(P)\Delta$ and
\begin{equation}\label{eq:example-fit}
 r_B(P)=\sum_{i=1}^4 h_i(P)
       =\sum_{i=1}^3 h_i(P)-\min_{p\in P}\sum_{i=1}^3(Bp)_i .
\end{equation}
This is the tetrahedral case of the simplex containment formula
of Nevskii~\cite{Nev}. We write $r(P)$ when $B$ is fixed. Translating $P$ by $u$ adds $g_i(u)$ to its
$i$th maximum. Since the four functions sum to zero, $r_B(P+u)=r_B(P)$.
For example, at $B=I$ every pair from $\{0,e_1,e_2,e_3\}$ has radius
one, whereas all four points together have radius three.

The covering radius $\mu=\mu(\Delta,L)$ is the least $t$ for which
$t\Delta+L=\R^3$. At $t=\mu$, Theorem~\ref{thm:closed-chain} gives
a degree-one cycle whose terms all have radius at most $\mu$.
Conversely, let $z$ be any periodic lattice $3$-cycle of degree one
and let $t$ be the largest radius of its terms. Every term is
admissible for $t\Delta$, so the same theorem gives $t\Delta+L=\R^3$
and hence $t\ge\mu$. Thus
\begin{equation}\label{eq:bottleneck}
 \mu(\Delta,L)=\min_{\partial z=0,\ \deg z=1}
                       \max_{P\in\operatorname{supp}z} r_B(P),
\end{equation}
where $P$ runs through the vertex sets of the terms with nonzero coefficient in $z$.
Degrees here refer to the actual points $Bp$. By
\eqref{eq:cycle-degree}, the normalized determinant of a term is its
integer determinant multiplied by $\operatorname{sgn}(\det B)$.

Call $B$ \reviewterm{generic} if no nonzero vector of $L$ is parallel to
a facet of $\Delta$. Equivalently, $g_i(u)\ne0$ for every
$u\in\Z^3\setminus\{0\}$ and each $i$. Two lattice points in a plane
parallel to facet $i$ would have a difference parallel to that facet.
Thus every such plane contains at most one lattice point.

\begin{lemma}\label{lem:perturb}
If $\Delta+B\Z^3$ covers with density less than $\alpha$, there are
generic covering bases $B_j\to B$ with density less than $\alpha$.
\end{lemma}
\begin{proof}
Take a degree-one cycle $z$ with $r_B(P)\le1$ on its support.
Choose $s<1$ so close to one that $1/(6s^3|\det B|)<\alpha$.
At $sB$ every term has radius $r_{sB}(P)=s\,r_B(P)\le s<1$.
The finitely many radii and the determinant vary continuously.
Near $sB$ the determinant keeps its sign, so the same chain of
integer labels still has degree one. Its radii remain below one,
and the covering density stays below $\alpha$.
For fixed $i$ and $u\in\Z^3\setminus\{0\}$, the equation $g_i(u)=0$
defines a proper hyperplane of matrices. Their countable union has
measure zero, so we can choose generic bases arbitrarily close to $sB$
and let $s$ tend to one.
\end{proof}

\subsection{The contact complex and its application to covering}
Fix a generic basis $B$. The complex of B\'ar\'any, Howe, and
Scarf~\cite{BHS} gives a degree-one cycle attaining the minimum
in~\eqref{eq:bottleneck}. To describe its cells, consider a translate
of a positive dilate of $-\Delta$ that is \reviewterm{maximal lattice-free}.
This means that its interior contains
no lattice point and that every strictly larger closed convex body
containing it has a lattice point in its interior. Write $x=By$,
so that $y$ gives the coordinates in the lattice basis. With $b_i$
the bound on the facet functional $g_i$, the tetrahedron is $E(b)$
in the $y$-coordinates and $BE(b)$ in the $x$-coordinates, with
\[
 E(b)=\{y\in\R^3:g_i(y)\le b_i,\ 1\le i\le4\},\qquad
 BE(b)=(b_1,b_2,b_3)-\Bigl(\sum_i b_i\Bigr)\Delta ,
\]
where $\sum_i b_i>0$. Consider one facet of $E(b)$.
If its relative interior contains no integer point, then every integer
point on it also lies on another facet. Moving just this facet by a
sufficiently small distance keeps those points on the boundary and
admits no new integer points, since only finitely many lie in the
bounded tetrahedron obtained by first increasing its bound by one.
This contradicts maximality.

Thus maximality gives at least one integer point in the relative
interior of each facet. Genericity gives at most one on its entire
plane. Write the points in \reviewterm{facet order} as $w_1,\ldots,w_4$,
so $w_i$ lies on facet $i$ of $E(b)$.
We call the lattice point $Bw_i$ in the relative interior of facet $i$
of $BE(b)$ its \reviewterm{facet contact}, and $w_i$ its contact label.
Thus the tuple $(w_1,w_2,w_3,w_4)$ specifies a vertex label for each facet. In particular
\[
                   g_i(w_j)<g_i(w_i)=b_i\qquad(j\ne i).
\]
Conversely, suppose we are given the set of four contact labels
$P=\{w_1,w_2,w_3,w_4\}$. Since $B$ is fixed, the four facet directions
are known. To locate the first facet in standard coordinates, observe
that every point of $BP$ satisfies $x_1\le b_1$, with equality at the
contact on that facet. Thus $b_1=\max_{w\in P}(Bw)_1$. The same argument
in each facet direction gives
\[
       b_i=\max_{w\in P}g_i(w)=h_i(P)\qquad (1\le i\le4).
\]
Taking these maxima requires only the set $P$ and also recovers the facet
assignment, since the maximizer for each $i$ is unique. The intersection
of the resulting half-spaces is the \reviewterm{container} of $P$.
Using $r_B(P)=\sum_i h_i(P)$ from~\eqref{eq:example-fit}, we obtain
\[
 E(b)=\bigcap_{i=1}^4\{y:g_i(y)\le h_i(P)\},\qquad
 BE(b)=(h_1(P),h_2(P),h_3(P))-r_B(P)\Delta.
\]

Use the four labels as the vertices of an abstract tetrahedron,
called a \reviewterm{cell}. Glue the cells along faces with the same vertex
labels to form the periodic \reviewterm{contact complex} $\scarf$.
The lattice translates of a cell form its \reviewterm{cell orbit}.
A cell is \reviewterm{solid}
if its four contact points are affinely independent and \reviewterm{flat}
if they are coplanar. The container $E(b)$ is three-dimensional in
both cases. In the flat case each point is exposed by its
own facet functional, so the four points are the vertices of a convex
quadrilateral. Its two diagonal subdivisions give four triangles in
all. They are the images of the four abstract faces obtained by omitting
one label at a time. The cell has zero volume, but each of these four
faces occurs in its boundary.

Let $A$ have rows $g_1,\ldots,g_4$. It has rank three and its rows
sum to zero, so $(1,1,1,1)$ spans its left kernel. Moreover,
$g_i(u)\ne0$ for each nonzero $u\in\Z^3$ by genericity.
These are the two hypotheses of~\cite{BHS}. Since $E(b)=\{y:Ay\le b\}$,
the contact labels defining our cells are exactly the vertices of the
maximal simplices of $K(A)$ in~\cite[p.~274]{BHS}. Thus $\scarf=K(A)$.
The map $w\mapsto Bw$, extended affinely over each cell, takes translation
by $u$ to translation by $Bu$ and induces $a_B:|\scarf|/\Z^3\to T_L$.

\begin{lemma}\label{lem:topology}
For a generic basis $B$, the contact cells have the following properties.
\begin{enumerate}
\item $E(b)\cap\Z^3=\{w_1,w_2,w_3,w_4\}$ for each maximal
lattice-free $E(b)$.
\item Every triangle belongs to two cells. Orienting each cell by listing
its vertices in facet order gives opposite induced orientations on their
common triangle. Around any fixed edge, the cells form one cycle,
consecutive cells sharing a triangle containing the edge.
\item In $\scarf/\Z^3$, sum the coherently oriented cells with coefficient
one, taking each translation class once. The resulting chain is a fundamental
cycle and has degree $\pm1$ under $a_B$.
\end{enumerate}
\end{lemma}

For (1), the interior of $E(b)$ is lattice-free and each facet plane
contains only its contact label~\cite[p.~274]{BHS}. In (2), the two
opposite incidences are Scarf's result recalled on the same page.
In the convex realization of~\cite[Theorem~3]{BHS}, the
\reviewterm{link} of an edge of $\scarf$ is the boundary of a
polygon whose sides correspond to the incident tetrahedra. Their cyclic
order gives the second assertion of (2).

For (3), $|\scarf|$ is homeomorphic to $\R^3$ and its quotient is a
three-torus with finitely many cells~\cite[Theorems~1--2 and Lemma~4]{BHS}.
Summing its coherently oriented top-dimensional cells gives the
fundamental cycle of the quotient, which we call the \reviewterm{contact cycle}.
In lattice coordinates their signed volumes sum to $\pm1$ by
P.~White's formula~\cite[p.~275]{BHS}. Applying $B$ multiplies this sum
by $\det B$, and the degree formula~\eqref{eq:cycle-degree} divides it by
$\covol(L)=|\det B|$.

\begin{corollary}\label{cor:scarf-radius}
For a generic basis $B$, the covering radius of $\Delta$ is the maximum
of $r(P)$ over the finitely many translation classes of cells of
$\scarf$. The fundamental cycle of $\scarf/\Z^3$, oriented to have
degree one under $a_B$, attains \eqref{eq:bottleneck}.
\end{corollary}
\begin{proof}
For each cell $P$, the tetrahedron $BE(b)$ is a translate of
$-r(P)\Delta$ with lattice-free interior. If $0<t<r(P)$, contraction
about an interior point places a closed translate $x-t\Delta$ in this
interior. Thus $(x-t\Delta)\cap L=\varnothing$, or equivalently
$x\notin t\Delta+L$. It follows that $\mu(\Delta,L)\ge\max_P r(P)$.

Conversely, take $t=\max_P r(P)$ and let $N$ be the admissible
complex for $t\Delta$ and $L$. For each cell $P$ of $\scarf$,
the lattice points $Bw$, $w\in P$, span a simplex of $N$ since $r(P)\le t$.
Thus $w\mapsto Bw$ identifies $\scarf$ with a subcomplex of $N$,
carrying translation by $u\in\Z^3$ to translation by $Bu\in L$.
On the quotients, the map $a:|N|/L\to T_L$ restricts to $a_B$.
Hence, by Lemma~\ref{lem:topology}(3), the image of the
fundamental cycle of $\scarf/\Z^3$, suitably oriented, is an
admissible cycle of degree one in $|N|/L$.
Theorem~\ref{thm:closed-chain} now gives $t\Delta+L=\R^3$, and this
cycle attains the minimum in~\eqref{eq:bottleneck}.
\end{proof}

To follow this cycle across a triangle with label set $F$, write
$F\cup\{d\}$ and $F\cup\{d'\}$ for the two cells incident to $F$
(Lemma~\ref{lem:topology}(2)). Their distinct fourth labels $d$ and
$d'$ need not lie on opposite sides of the plane of $F$.
Scarf states the adjacency rule using row minima for inequalities
$Ay\ge b$~\cite[Section~3]{Scarf}. With our inequalities
$g_i(y)\le b_i$, the contact on each facet gives a row maximum.

\begin{lemma}[Scarf's adjacency rule in dimension three]\label{lem:transfer}
Suppose $d$ is the contact label on facet $i$ of the first container.
Let $w$ be the unique maximizer of $g_i$ on $F$, and let $j$ be its
facet number in the first container. Then $j\ne i$. In the second
container $w$ is the contact label on facet $i$ and $d'$ is the
contact label on facet $j$. The other two labels in $F$ keep their
facet numbers.
\end{lemma}

This transfer also explains the orientation in Lemma~\ref{lem:topology}(2).
Replacing $d$ by $d'$ in slot $i$ gives the same induced orientation on
$F$. Putting the new cell in facet order then interchanges slots $i$ and
$j$, reversing that orientation.
The two cells on $F$ represent different translation classes. A
translation between them would have to fix each of the two contact
labels that keep their facet numbers, hence would be zero.
When the basis is unknown, the three labels in $F$ give at most three
choices for the facet number $j$ of $d'$. To locate such neighbours, we
return to the question in the introduction. In a generic covering of density
less than two, must some solid cell have three edges at one vertex that form
a basis of $L$? Translating that vertex to the origin would then make the
cell the coordinate simplex in this basis.

\section{Volume and index restrictions on the cells}\label{sec:highindex}

A covering of density less than $125/63$ would, by
Lemma~\ref{lem:perturb}, give one in general position. For the volume
bounds we consider a generic covering $\Delta+L$ of density $\theta<2$,
where $L=B\Z^3$.

\subsection{Cell indices and volume bounds}
Let $P=\{p_0,p_1,p_2,p_3\}$ be the label set of a solid cell, so
its four contact points are $Bp_0,\ldots,Bp_3$ in $L$. Its three edges
from $Bp_0$ form a basis of $L$ precisely when the differences $p_i-p_0$,
$1\le i\le3$, generate $\Z^3$. The \reviewterm{index} of the cell is
the index of the sublattice generated by these differences in $\Z^3$.
It equals the \reviewterm{normalized volume} of the cell in these coordinates,
\[
 \idx(P)=|\det[p_1-p_0,p_2-p_0,p_3-p_0]|=6\vol(\conv P).
\]
Thus an index-one cell gives the desired basis. In the original
coordinates the same index is $6\vol(\conv BP)/\covol(L)$.

In the lattice basis $B$, the reflected tetrahedron $-\Delta$ becomes
\[
          K_B=-B^{-1}\Delta,\qquad
          \vol(K_B)=\frac{1}{6|\det B|}=\theta .
\]
Since $\Delta+L$ covers, $K_B+\Z^3$ covers as well. By
Corollary~\ref{cor:scarf-radius}, the four vertices of each contact
cell $P$ can be translated into $K_B$. The translation may depend on
$P$, but the tetrahedron $K_B$ is the same for every cell.

Let $p,q\in\Z^3$ be two vertices of the same contact cell. They lie
in a common translate of $K_B$, so $q-p\in K_B-K_B$. Call a lattice
vector \reviewterm{short} if it belongs to $K_B-K_B$. Since the same $K_B$
serves every cell, we can bound the determinant of three edges even
when they belong to different cells. For any tetrahedron $T$ of
positive volume,
\begin{equation}\label{eq:rootdet}
 \max_{u,v,w\in T-T}|\det[u,v,w]|=6\vol(T).
\end{equation}
If $T=x+A\Delta$, then $T-T=A(\Delta-\Delta)$, and $A$ multiplies
both absolute determinants and volumes by $|\det A|$. It therefore
suffices to prove that the maximum is one for $\Delta$.
Put $e_0=0$. The vertices of $\Delta-\Delta$ are the
differences $e_i-e_j$, $0\le i\ne j\le3$. Regard each such vector as
an arrow from $j$ to $i$ on four points labelled $0,1,2,3$. Three
such vectors give the columns of a determinant. If two underlying
edges coincide or the three form a triangle, the columns are
dependent. Otherwise the three edges form a tree on the four points.
At an endpoint other than $0$, the corresponding row has a single
nonzero entry, equal to $1$ or $-1$. Expanding along this row removes
the endpoint and its edge. Repetition gives absolute determinant one.
Multilinearity extends the bound to the whole difference body, and
$e_1,e_2,e_3$ give equality. Taking $T=K_B$ in~\eqref{eq:rootdet},
we obtain $|\det[u,v,w]|\le6\vol(K_B)=6\theta<12$ for any three
short lattice vectors. Their determinant is an integer, so its
absolute value is at most $11$. Taking the three edges at one vertex
of a solid cell gives $\idx(P)\le11$.

Can three cells occur in the same covering if each is small enough
on its own? In the example below their normalized volumes are only
$6,1,1$, yet together they force $\theta\ge5/2$. The extra information
comes from having independent translates of all three hulls inside the
same tetrahedron $K_B$. To use it, let $H_1,H_2,H_3$ be their hulls.
The volume of their Minkowski sum $sH_1+tH_2+uH_3$, for $s,t,u\ge0$,
is a homogeneous cubic polynomial. Its coefficient of $stu$ is
$6V(H_1,H_2,H_3)$, where $V$ is the \reviewterm{mixed volume}. In particular
$V(H,H,H)=\vol(H)$. Translating the three hulls independently
translates their sum and hence preserves every coefficient. The mixed
volume also increases weakly when any of its arguments is enlarged
\cite[Proposition~2.1]{ABS}.

\begin{lemma}[Mixed volumes]\label{lem:mix}
For three cells $P_1,P_2,P_3$ of a generic covering of density
$\theta<2$, put $H_i=\conv P_i$ in integer coordinates. Then
$6V(H_1,H_2,H_3)\le11$. The cells may repeat or be flat.
\end{lemma}
\begin{proof}
Choose $t_i$ with $H_i+t_i\subset K_B$. Independent translation and
then enlargement give
\[
 V(H_1,H_2,H_3)
 =V(H_1+t_1,H_2+t_2,H_3+t_3)
 \le V(K_B,K_B,K_B)=\theta .
\]
Since the $H_i$ have integer vertices, $6V(H_1,H_2,H_3)$ is an
integer, also for lower-dimensional $H_i$ \cite[Proposition~2.1(4)]{ABS}.
It is at most $6\theta<12$.
\end{proof}

To see how three hulls give a stronger bound, take
\[
\begin{aligned}
 H_1&=\conv\{0,(1,0,0),(0,0,1),(1,6,1)\},\\
 H_2&=\conv\{0,(0,-1,0),(1,4,0),(1,6,1)\},\\
 H_3&=\conv\{0,(0,1,-1),(1,4,0),(0,-1,0)\}.
\end{aligned}
\]
Their individual normalized volumes $6\vol(H_i)$ are $6,1,1$, each
within the bound. The normalized volumes of $H_1+H_2$, $H_1+H_3$,
$H_2+H_3$ and $H_1+H_2+H_3$ are $58,67,32,239$. Taking the alternating
sum of the ordinary volumes gives
\[
 6V(H_1,H_2,H_3)=(239-58-67-32+6+1+1)/6=15.
\]
Thus the three hulls together force $\theta\ge 5/2$.

\subsection{How large can a cell index be?}
The shape of a single solid cell gives a first improvement on the bound
$\idx(P)\le11$. To sharpen it further, we ask whether the cell can occur
with all the neighbours required by the contact complex.

\begin{proposition}\label{prop:high}
In a generic tetrahedron lattice covering of density less than $2$, every
solid cell has index at most $8$. If all cells are solid, every index is
at most $5$.
\end{proposition}
\begin{proof}
\emph{The shape of one cell.}
Let $H$ be the convex hull of a solid cell of index $m>1$.
By Lemma~\ref{lem:topology}(1), its only lattice points are its vertices.
By G.~K.~White's theorem~\cite[Theorem~1]{White}, two opposite edges of $H$
lie in consecutive lattice planes. Take one of these primitive edges
as $e_1$, complete it to a basis $e_1,e_2$ of the plane lattice, and
take a joining edge from the same vertex as $e_3$. With this vertex
at the origin and a choice of sign for $e_2$, we obtain
\[
 H=H_{m,t}=\conv\{0,e_1,e_3,(t,m,1)\},\qquad \gcd(t,m)=1 .
\]
The last vertex has height one, its second coordinate is $m$ by the
index formula, and the upper edge $(t,m,0)$ is primitive. We use
this same lattice basis for $K_B$ and all the other cells.

\smallskip\noindent
\emph{The density forced by this cell.}
We know that a translate of $K_B$ contains $H$. In the other
direction, let $r$ be the least factor for which a translate of
$rH$ contains $K_B$. We shall show that
\[
                   1+\frac{2}{m}\le r\le\frac{6\theta}{m}.
\]
The upper bound follows from mixed volume and the lower bound from
the covering by $rH+\Z^3$.

For both estimates let $W$ have columns $e_1,e_3,(t,m,1)$, so
$H=W\Delta$ and $|\det W|=m$. For any polytope $S$ set $Y=W^{-1}S$.
Write $r_H(S)$ for the least $t\ge0$ such that some translate of $tH$
contains $S$. In the coordinates of $Y$ this gives
\[
 r_H(S)=\max_{y\in Y}\sum_i y_i-\sum_i\min_{y\in Y} y_i,
 \qquad 6V(H,H,S)=m\,r_H(S).
\]
Indeed, the lower coordinate planes of a smallest translate of
$r_H(S)\Delta$ around $Y$ pass through its coordinate minima. The
opposite plane must reach the maximum of $\sum_i y_i$, which gives
the first formula.

For the second formula, translate $Y$ so that its minimum in each coordinate
is zero. This preserves both $r_H(S)$ and $V(\Delta,\Delta,Y)$, and now
$r_H(S)=\max_{y\in Y}(y_1+y_2+y_3)$. The facet formula for
$3V(\Delta,\Delta,Y)$ sums, over the facets of $\Delta$, the area of each
facet times the maximum of $u\cdot y$ on $Y$ for its unit outer normal
$u$~\cite[eq.~(2.8)]{ABS}. For the three coordinate facets these maxima
are zero. The remaining facet has area $\sqrt{3}/2$ and unit outer normal
$(1,1,1)/\sqrt{3}$, so
\[
  3V(\Delta,\Delta,Y)
       =\frac{\sqrt{3}}{2}\,\frac{r_H(S)}{\sqrt{3}}=\frac{r_H(S)}{2}.
\]
Applying $W$ multiplies mixed volume by $|\det W|=m$. Hence
$6V(H,H,S)=m\,r_H(S)$.

With $S=K_B$ we have $r=r_H(K_B)$, and two copies of $H$ can be
translated into $K_B$. Consequently $mr=6V(H,H,K_B)\le6\theta$.
Since $K_B+\Z^3$ covers and $K_B$ fits in a translate of $rH$,
$r\Delta+\Lambda$ covers for $\Lambda=W^{-1}\Z^3$. Every point
of $\Lambda$ has coordinates in $m^{-1}\Z$, and its last two
coordinates add to an integer, the original third coordinate.
Every point $y$ in the open box
\[
                    Q=(0,1/m)\times(0,1)\times(0,1/m)
\]
must lie in $\lambda+r\Delta$ for some $\lambda\in\Lambda$.
Then $\lambda_i\le y_i$ for all $i$ and
$\sum_i(y_i-\lambda_i)\le r$. Since $y_1,y_3<1/m$ and
$\lambda_1,\lambda_3\in m^{-1}\Z$, we have $\lambda_1,\lambda_3\le0$.
Also $\lambda_2\le y_2<1$, so the integer $\lambda_2+\lambda_3$ is
at most zero. Hence $\sum_i\lambda_i\le0$ and $r\ge\sum_i y_i$
throughout $Q$. Taking the supremum over the box gives
\[
       r\ge\sup_{y\in Q}(y_1+y_2+y_3)=1+2/m.
\]
Together with $mr\le6\theta$ this gives $m+2\le6\theta<12$.
Already the single cell forces $m\le9$.

For $6\le m\le9$, which of the parameters $t$ coprime to $m$ must we test?
In the plane lattice $M=\Z e_1+\Z e_2$, the lower and upper edge vectors
are $u=e_1$ and $v=tu+me_2$. Replacing $e_2$ by $e_2+ku$ for $k\in\Z$
writes the same upper edge with coefficient $t-km$, so we take $t$ modulo $m$.
Reversing the upper edge changes $t$ to $-t$. We may also interchange
the two planes and take the primitive vector $v$ as the first basis
vector of $M$. The original relation $v\equiv tu\pmod{mM}$ then reads
$u\equiv t^{-1}v\pmod{mM}$, so the new parameter is $t^{-1}$ modulo $m$.
These choices give different descriptions of the same cell.

For $m=7$ the class of $1$ is $\{1,6\}$. Since $2\cdot4\equiv1\pmod 7$,
the residues $2,4$ and their negatives $5,3$ form the other class.

\begin{samepage}
Taking one representative per class for each $m$ gives
\[
\begin{array}{c|cccc}
 m&6&7&8&9\\ \hline
 t&1&1,2&1,3&1,2
\end{array}
\]
\end{samepage}

\smallskip\noindent
\emph{Which neighbours can occur?}
Fix one of these pairs $(m,t)$ and suppose that the covering contains a
cell with hull $H=H_{m,t}$. Keep the lattice basis chosen for $H$.
For every other cell, let $P\subset\Z^3$ be its set of four vertex
labels in this basis. This records both the shape of its hull and its
orientation relative to the coordinate axes.
Apply the mixed-volume formula above with $S=\conv P$. Writing
$a(v)=mW^{-1}v=(mv_1-tv_2,mv_3-v_2,v_2)$, Lemma~\ref{lem:mix} gives
\begin{equation}\label{eq:highwindow}
 6V(H,H,\conv P)=\max_{v\in P}{\textstyle\sum_i}a_i(v)
                       -\sum_i\min_{v\in P}a_i(v)\le11 .
\end{equation}
Subtract $\min_{v\in P}a_i(v)$ in each coordinate. The three resulting
coordinates are nonnegative, and their sum at each vertex is at most
eleven by~\eqref{eq:highwindow}. Thus each coordinate has range at
most eleven. After translating one vertex to zero, every vertex satisfies
$a(v)\in[-11,11]^3$. Enumerating the vertex sets in this bound gives a
finite list of candidate cells up to translation for this $(m,t)$.

Can we select cells from this list, including $H$, that form the contact
cycle of a covering? Initially only the cell with hull $H$ is prescribed.
A triangular face that occurs just once among the prescribed cells,
counting modulo translation, is an \reviewterm{unpaired face}. Its second
incidence must come from a neighbour. To extend a choice, take such a face
and try each candidate that shares it up to translation.
Lemma~\ref{lem:transfer} restricts the possible facet assignments.
All chosen assignments must be compatible with the same four facet
functionals, and every triple of chosen cells must satisfy Lemma~\ref{lem:mix}.
Each choice of neighbour gives a branch in which that cell is also prescribed.

At every stage we can test two necessary conditions on a completion.
First, the completed chain must be closed, so the number of occurrences of
each triangle class in its boundary sum is even. Second, by
Lemma~\ref{lem:topology}(3) the contact cycle has degree $\pm1$ under
$a_B$. Recall that to define this map we send each vertex label $w$ to
$Bw$, extend affinely over each cell, and project to the torus $T_L$.
Reversing the orientation of the cycle if necessary, take its degree to
be one. Choose a point $y$ of the torus outside the images of
all candidate faces and flat cells. Every preimage of $y$ then lies in
the interior of a solid cell. If $N_+$ and $N_-$ count preimages where
$a_B$ preserves and reverses orientation, respectively, its degree gives
$N_+-N_-=1$. Thus the total $N_++N_-=1+2N_-$ is odd.

For each candidate $P$ let $x_P$ be $1$ if it is selected and $0$
otherwise. Fix $x_P=1$ for each cell already prescribed on the branch.
To compute the degree modulo two from these indicators, order the four
vertex labels of each cell lexicographically in the fixed basis. Compare
the first coordinates, then the second in case of a tie, and
finally the third. For $H_{7,2}$ the resulting order is
\[
 p_0=(0,0,0)<p_1=(0,0,1)<p_2=(1,0,0)<p_3=(2,7,1).
\]
The subscripts give the positions in this ordering independently of the
facet assignments. The order agrees on common faces and is preserved
by translation.
For $j=1,2,3$, assign to each oriented edge $[u,v]$ the number
$\alpha_j([u,v])=(v_j-u_j)\bmod 2$. Such an assignment of a value to
each edge is an \reviewterm{edge cochain}. Its value is one exactly when
the $j$th coordinate changes by an odd integer.

\begin{samepage}
With the vertices of $P$ in lexicographic order, the \reviewterm{simplicial cup product}
of three edge cochains multiplies their values on the successive edges
$[p_0,p_1]$, $[p_1,p_2]$ and $[p_2,p_3]$, respectively.
For the three coordinate cochains this gives the coefficient
\[
 \begin{aligned}
 c(P)&=(\alpha_1\smile\alpha_2\smile\alpha_3)(P)\\
     &=(p_1-p_0)_1(p_2-p_1)_2(p_3-p_2)_3\bmod 2.
 \end{aligned}
\]
\end{samepage}
To see why this weight counts degree, use $B^{-1}$ to identify $T_L$
with $\R^3/\Z^3$. For a closed edge path, lift its image under $a_B$
to a continuous path in $\R^3$ by starting at one representative of its
first vertex and following the edge increments $v-u$. The resulting
path is a \reviewterm{lift} because its projection modulo $\Z^3$ is
the given torus path. Its endpoints differ by an integer vector whose $j$th
coordinate is the \reviewterm{winding number} around the $j$th
coordinate circle. Thus $\alpha_j$ sums to this winding number modulo
two on any closed path. It sums to zero around each triangle, so it is
a cocycle representing the pullback of the $j$th coordinate class.

The product of the three coordinate classes evaluates to one on the
fundamental cycle of the torus. Indeed, among the six Kuhn tetrahedra
of the unit cube, only the one with successive steps $e_1,e_2,e_3$
contributes to this product. Cup products commute with pullback
by~\cite[Proposition~3.10]{Hatcher}.
Consequently, on a closed selection the sum $\sum_P c(P)x_P$ evaluates
the product on the image cycle and gives its degree modulo two. This is
also the parity of the preimage count above. Indeed, let
$\widetilde y\in\R^3$ represent $B^{-1}y$. Each solid $P$ contributes one
preimage for each point of $(\widetilde y+\Z^3)\cap\operatorname{int}(\conv P)$,
and each flat $P$ contributes none. For the contact cycle this parity is one.
Together with face cancellation this gives
\begin{equation}\label{eq:cup}
 \sum_{P\supset F}x_P\equiv0\pmod2,\qquad
 \sum_P c(P)x_P\equiv1\pmod2,\qquad x_P=1\quad(P\text{ prescribed}).
\end{equation}
The first congruence counts the occurrences of each triangle class $F$
with multiplicity. If the system has no solution, no covering can contain
all the cells prescribed on this branch.

The search in Appendix~\ref{sec:high-replay} excludes both representatives
at $m=9$ even when flat cells are allowed. In the all-solid case it also
excludes $m=6,7,8$. Together with $m\le9$, this proves the two bounds.
\end{proof}

\subsection{The greatest common divisor of the cell indices}
The upper bound alone would still allow all solid indices to have a
common prime factor. Adjacency rules this out for generic coverings of
any density, even in the presence of flat cells.

\begin{lemma}\label{lem:analytic-index-gcd}
The indices of the solid cells of a generic tetrahedron lattice
covering have greatest common divisor one.
\end{lemma}
\begin{proof}
At least one cell is solid, since the oriented cells have nonzero
total signed volume. Take any face with vertices $a,b,c$ and let
$d$ be the fourth contact of one of its cells. The triangle $abc$
lies in the closed lattice-free tetrahedron of these four contacts.
By Lemma~\ref{lem:topology}(1) the only possible extra lattice
point on the triangle is $d$. The facet functional of $d$ is
strictly larger at $d$ than at $a,b,c$, so $d\notin\conv\{a,b,c\}$.
This also proves emptiness of the face when the four contacts are
coplanar.

In the lattice of integer vectors parallel to the face, Pick's formula
gives the area $0+3/2-1=1/2$, since the triangle's only lattice points
are its vertices. Thus $b-a,c-a$ form a basis of this plane lattice.
A primitive integer normal
to the plane takes the value one on some integer vector. Adding
this vector extends the two edges to a basis of $\Z^3$. In these
coordinates, with $a$ as origin, the face is $0,e_1,e_2$ and the
fourth vertex has the form $(u,v,h)$. The cell determinant is $h$
up to sign. In particular, for any prime $p$ the common face has
three noncollinear vertices modulo $p$, and if $p\mid h$ the
fourth vertex belongs to their plane $x_3=0$ modulo $p$.

Suppose $p$ divides every solid index. Each flat determinant is
already zero. Thus in every cell the residues of its four labels
lie in the unique affine plane determined by any three on a face.
Two cells sharing a face give the same plane. The abstract complex before quotienting by translations triangulates a
connected $3$-manifold. A path between tetrahedron interiors can be taken
to cross triangular faces away from the edges and vertices. Thus any two
cells are joined by a sequence of cells sharing faces. Hence all vertex
labels lie in one affine plane modulo $p$. Such a plane contains only
$p^2$ residue classes, whereas the integer translates of any one vertex
represent all $p^3$ residue classes modulo $p$. This gives a contradiction.
\end{proof}

When all cells are solid and $\theta<2$, every index is at most five.
If the least index were five, all indices would equal five, contrary to
the lemma. Section~\ref{sec:noflat} rules out least indices two, three
and four. When flat cells occur,
Section~\ref{sec:flat-components} finds a lattice basis at a triangle
shared by a flat cell and a solid cell. In either basis the volume
bounds restrict the possible positions of the remaining contacts.
Section~\ref{sec:finite} then imposes the containment inequalities
for all their cells on the same basis matrix $B$.

\section{Coverings with only solid cells}\label{sec:noflat}
Suppose all cells are solid. By Lemma~\ref{lem:analytic-index-gcd}
their indices have greatest common divisor one. The pair of indices two
and three already has this property. To find a cell of index one we use
the restrictions on which cells can share a triangular face. The three
edges at a vertex of such a cell then give a lattice basis in which
we can describe the whole complex.

For a cell with vertex set $P=\{w_1,w_2,w_3,w_4\}$, record its contact
labels in facet order as $(w_1,w_2,w_3,w_4)$, so that $Bw_i$ lies on
facet $i$ of the container in standard coordinates. To compute
$c(P)$ in~\eqref{eq:cup}, we list these same vertices in lexicographic order.
For a chain $z$, $\operatorname{supp}z$ consists of the terms with nonzero
coefficient, as in~\eqref{eq:bottleneck}. In the classification below, a
\reviewterm{support} records the cell orbits together with the facet assignment
on each cell. We can translate
each tuple independently and list the cells in any order. A change of
lattice basis or a permutation of the facets acts on all the cells at once.
Two supports can give the same degree-one cycle. Their facet assignments
specify which vertex $w_i$ must satisfy $g_i(w_i)>g_i(v)$ for every
other vertex $v$ of its cell. These comparisons restrict the possible
neighbours in Lemma~\ref{lem:transfer} and the basis matrices $B$ in
Section~\ref{sec:finite}.

\begin{proposition}\label{prop:solid-catalogue}
Suppose all maximal cells of a generic tetrahedron lattice covering
of density $\theta<2$ are solid. Some cell has index one. Up to a
change of lattice basis and a common permutation of the facets, the
cells with their facet assignments form one of the three supports displayed
in Section~\ref{sec:solid-supports}. They have six, six and fourteen cell
orbits, respectively.
\end{proposition}

We prove Proposition~\ref{prop:solid-catalogue} in two steps.
Section~\ref{sec:cell-index-one} rules out each possible minimum cell
index greater than one. Section~\ref{sec:solid-supports} then takes an
index-one cell as the coordinate simplex and pairs triangles to determine
the full support.

\subsection{A cell of index one}\label{sec:cell-index-one}
Work with the covering $\Delta+L$ and put $D=\covol(L)$. Write $C\subset L$
for the four lattice points of a cell, and let $q$ be the least cell index.
Under $\theta<2$, Proposition~\ref{prop:high} gives
\begin{equation}\label{eq:noflat-indices}
  q\le \idx(C)=6\vol(\conv C)/D\le5
  \qquad\text{for every maximal cell }C .
\end{equation}
The indices have greatest common divisor one by
Lemma~\ref{lem:analytic-index-gcd}. Hence $q=5$ is impossible.
To rule out $q=2,3,4$, we determine the lattice in coordinates of a
minimum cell and then pair its triangles. For $q=3,4$ we can forget
the facet assignments and test the signed volumes when all triangles are
paired. For $q=2$ we also use which vertex meets each facet.

Translate one vertex of a cell of index $q$ to zero and write the
other three points of $L$ as the columns of $W$. The group
$W^{-1}L/\Z^3$ has order $q$. Each coordinate
projection to $\R/\Z$ is injective. Indeed, if a nonzero class has
first coordinate zero, then it or its negative has a representative
$(0,s,t)$ with $s,t\ge0$ and $s+t\le1$. This would give a further
lattice point on a face. Changing the vertex at the origin gives the
same conclusion for all four barycentric coordinates. In particular
the group is cyclic. For $q>1$ choose its generator to have residue
$1$ in the first Cartesian coordinate. Its four barycentric residues,
in units of $1/q$, are
\[
                      (-1-a-b,1,a,b)\pmod q .
\]
For every nonzero multiple the positive residues sum to $2q$.
Indeed the sum is a multiple of $q$, and a sum of $q$ or $3q$
puts the class or its negative in the closed simplex. The possible pairs $(a,b)$ are
\[
\begin{array}{c|l}
 q& (a,b)\\ \hline
 2&(1,1)\\
 3&(1,2),(2,1),(2,2)\\
 4&(1,3),(3,1),(3,3).
\end{array}
\]
In each case a permutation of the vertices gives the generator
$(-1,1,1,-1)$. For each point $x\in L$, multiply its edge coordinates
$W^{-1}x$ by $q$ to obtain the integer vector $p=qW^{-1}x$.
For $q=1,2,3,4$, the lattice and the vertex set of the chosen cell in these coordinates are
\[
 \Gamma_q=qW^{-1}L=q\Z^3+\Z(1,1,-1),\qquad
 V=\{0,qe_1,qe_2,qe_3\}.
\]
The label $p\in\Gamma_q$ represents the point $Wp/q\in L$. Since
$D=|\det W|/q$, the covolume of $\Gamma_q$ is $q^3D/|\det W|=q^2$. Number the facets so that the contact labels on $V$ are
$(0,qe_1,qe_2,qe_3)$ in facet order. A label $p$ is short when $Wp/q\in\Delta-\Delta$. Every cell
edge has a short label since its endpoints in $L$ lie in a common
translate of $-\Delta$.
For any three short labels $u,v,w\in\Gamma_q$, applying
\eqref{eq:rootdet} with $T=\Delta$ gives
\begin{equation}\label{eq:noflat-triples}
  \frac{|\det[u,v,w]|}{q^2}
  =\frac{|\det[Wu/q,Wv/q,Ww/q]|}{D}\le6\theta<12 .
\end{equation}
Let $p\in\Gamma_q$ be short. For each nonempty subset of $\{1,2,3\}$,
partition the vertices of $V$ into the indicated set and its
complement, with $0$ in the complement. Choose two edges of $V$ that together give a spanning tree on each
part, with a singleton requiring no edge. Their determinant with
$p$, divided by $q^2$, is the sum of the corresponding coordinates
of $p$, up to sign. For example,
\[
 |\det[p,q(e_1-e_2),qe_3]|/q^2=|p_1+p_2|.
\]
These sums are integers, so each has absolute value at most $11$.
Consequently every short label lies in the finite set
\[
                     U_q=\Gamma_q\cap11(\Delta-\Delta).
\]

At a triangle having only one prescribed incidence, translate a vertex
to zero and try every fourth vertex in $U_q$. Require all six edge
differences to lie in $U_q$, the index to lie in $[q,5]$, and the
closed hull to contain only its four lattice vertices. When the facet assignment on the initial cell
is prescribed, try each of the three transfers of Lemma~\ref{lem:transfer}
and impose the strict comparisons in both cells. These tests include
the neighbour of any actual covering extending the prescribed cells.
A connected selection in an actual covering with all triangles paired
is the whole quotient complex, since its dual graph (cells joined
across common triangles) is connected.
Orient its cells coherently. For each cell take the determinant of the three
edge vectors in that orientation and divide by $\covol(\Gamma_q)=q^2$.
The absolute value is the cell index. Their sum, the signed index sum, is six times the degree
in $\Gamma_q$-coordinates by~\eqref{eq:cycle-degree}, hence $\pm6$ by Lemma~\ref{lem:topology}(3).

Appendix~\ref{sec:minimum-indices} gives the two searches and their counts.
For $q=3,4$, indices, lattice emptiness, edge determinants and face
multiplicities suffice to continue from $V$ without prescribing
facet assignments. At $q=3$ there are four fully paired selections,
all with signed index sum zero. All other branches fail a local test,
including every branch at $q=4$. At $q=2$ the neighbour across a
fixed face gives $65$ possibilities up to symmetry before imposing
strict comparisons. The strict contact comparisons, lattice emptiness and
opposite incidences reject $63$ at the first step. The remaining two
are exhausted after adjoining at most two cells. This proves $q=1$.

\subsection{The three candidate supports}\label{sec:solid-supports}
We give the contact labels of each candidate in lattice coordinates.
For a support $S$ the corresponding periodic chain is
\[
           z(S)=\sum_{(w_1,w_2,w_3,w_4)\in S}[w_1,w_2,w_3,w_4],
\]
with one term per cell orbit, oriented in the order of the four facets.

\emph{Candidates I and II (six cells each).}
Put $s=e_1+e_2+e_3$. Both supports consist of the six Kuhn tetrahedra
of the unit cube, with vertex sets $\{0,e_i,e_i+e_j,s\}$ as $(i,j,k)$
runs through the permutations of $(1,2,3)$. Each entry in the table
gives the facet number at the vertex above it.
\[
\begin{array}{c|cccc}
   &0&e_i&e_i+e_j&s\\ \hline
 \mathrm{I}  &4&j&k&i\\
 \mathrm{II} &i&j&k&4
\end{array}
\]
For instance, on $\{0,e_1,e_1+e_2,s\}$ the tuple of contact labels in facet
order is $(s,e_1,e_1+e_2,0)$ in I and $(0,e_1,e_1+e_2,s)$ in II.
Thus interchanging the contacts at the ends of the body diagonal
reverses each oriented simplex. After orienting both chains to degree
one in these coordinates, we obtain the same cycle with two different
facet assignments.

\begin{samepage}
\emph{Candidate III (fourteen cells).}
Take the last cell as the coordinate simplex and translate each cell
so that its contact on facet one is $0$. Each row gives the remaining
contacts, so its contribution to
$z(S)$ is $[0,w_2,w_3,w_4]$.
\[
\begin{array}{r|ccc}
 \text{cell}&w_2&w_3&w_4\\ \hline
 1 & (-2,-1,1) & (-2,-1,2) & (-1,0,1)\\
 2 & (-2,-1,2) & (-2,0,1) & (-1,0,1)\\
 3 & (-1,-2,1) & (-1,-2,2) & (1,-1,0)\\
 4 & (-1,-2,2) & (-1,-1,1) & (1,-1,0)\\
 5 & (-1,-1,1) & (-2,-1,1) & (-1,0,1)\\
 6 & (-1,-1,1) & (-1,-1,2) & (0,1,0)\\
 7 & (-1,-1,1) & (-1,0,1) & (-1,-1,2)\\
 8 & (-1,-1,2) & (-1,0,1) & (0,1,0)\\
 9 & (0,-2,1) & (-1,-2,1) & (1,-1,0)\\
 10 & (0,-1,1) & (-1,-1,1) & (0,1,0)\\
 11 & (0,-1,1) & (0,1,0) & (1,0,0)\\
 12 & (1,-1,0) & (-1,-1,1) & (0,-1,1)\\
 13 & (1,0,0) & (0,0,1) & (1,1,0)\\
 14 & (1,0,0) & (0,1,0) & (0,0,1)
\end{array}
\]
\end{samepage}
The third cell has index three and the other thirteen have index one.
Modulo lattice translations, the triangular faces cancel in $12$ pairs
for each of I and II and in $28$ pairs for III. In the displayed
coordinates the signed determinants sum to $-6$, $6$ and $-6$,
respectively. Hence each chain has degree one up to its overall sign
by~\eqref{eq:cycle-degree}.

It remains to show that every covering satisfying the hypotheses of
Proposition~\ref{prop:solid-catalogue} has one of these three supports.
Choose an index-one cell as in the preceding subsection and take its
three edges at one vertex as a basis.
Translate one vertex of each unpaired face to zero. The fourth
vertex of its neighbour in an actual covering then belongs to $U_1$,
and this neighbour satisfies the tests of Appendix~\ref{sec:unit-runs}.
Thus some branch follows the covering until all triangles are paired.
At that point the selected cells form the whole quotient complex,
since its cells are connected across triangles. The enumeration in
Appendix~\ref{sec:unit-runs}, implemented in \texttt{solid/run.py},
yields twenty paired lists. Changing the
lattice basis and permuting the facets simultaneously on all cells
reduces them to I, II and III, proving Proposition~\ref{prop:solid-catalogue}.

\section{Components of coplanar cells}\label{sec:flat-components}
Suppose a generic covering has a flat cell. Consider the flat cells in
the fundamental $3$-cycle of the quotient contact complex
(Lemma~\ref{lem:topology}). Two flat cell orbits are adjacent if
representatives share a triangular face. We show that each connected
component of flat cells forms a \reviewterm{path} under this
adjacency. Its length is the number of cell orbits in the component.
When $\theta<2$, a solid neighbour at an end triangle supplies a lattice
basis. Section~\ref{sec:finite} starts with a longest such path and a
solid neighbour at its initial end, with their facet assignments.
Choosing a longest path bounds the length of every other flat component
in a completion.

\begin{proposition}\label{prop:flat-catalogue}
Suppose a generic tetrahedron lattice covering has a flat cell.
Modulo translations, its flat cells form paths of parallelograms under
triangle adjacency. If the density is $\theta<2$, each solid cell
meeting an end triangle of a path has index one, and a longest path
has one of the thirteen words in~\eqref{eq:flat-words}, with the
initial facet assignments specified below.
\end{proposition}

\subsection{Paths and their end triangles}
Write the covering as $T+L=\R^3$ and take the vertices in $L$.
Starting with a flat cell, follow its neighbours across triangular
faces as long as they too are flat. We shall show that the cells
reached in this way form a path and identify the triangles where the
path meets solid cells. Throughout we identify lattice translates.

\emph{The four faces of a flat cell.}
Each of its four contact points is exposed by a facet of its container,
so their convex hull is a quadrilateral. By Lemma~\ref{lem:topology}(1)
it contains no other lattice point. Translate one vertex to zero and
write the vertices in cyclic order as $0,a,d,b$. The lattice in this plane is $L\cap\operatorname{span}_{\R}\{a,b\}$. Each triple of vertices spans a triangle containing no other
lattice point. By Pick's formula its area
is half the covolume of the plane lattice, so its two edge vectors at a
vertex form a basis of this lattice. We call such
a triangle \reviewterm{unimodular}. In particular, $a,b$ are a basis
of the plane lattice. Compute planar determinants in this ordered basis,
so that $\det(a,b)=1$. Writing $d=ra+sb$, the two triangles on the
diagonal $[0,d]$ give $s=\det(a,d)=1$ and $r=\det(d,b)=1$. Hence the
cell has vertex set $P(a,b)=\{0,a,b,a+b\}$, the vertices of a
fundamental parallelogram of the plane lattice.

The four triangular faces of this cell are the triples obtained by
omitting one vertex. They occur in pairs. The diagonal $[a,b]$ gives
$S_1=\{0,a,b\}$ and $a+b-S_1=\{a,b,a+b\}$. The diagonal $[0,a+b]$
gives $S_2=\{0,a,a+b\}$ and $a+b-S_2=\{0,b,a+b\}$. Thus a flat cell
with a face in the translation class $[S]$ also has a face in $[-S]$.
These classes are distinct since a triangle has no centre of symmetry.

\emph{How flat cells are joined.}
For example, $P(a,b)$ and $P(a,a+b)=\{0,a,a+b,2a+b\}$ share the
face $S_2$. Both therefore also supply a face of class $[-S_2]$. To
record such adjacencies, fix the plane direction, which is determined
by any shared triangle. For each unimodular triangle $S$ use the pair
$N(S)=\{[S],[-S]\}$ as a \reviewterm{node} of a graph. Represent
$P(a,b)$ by the edge joining $N(S_1)$ to $N(S_2)$. The two drawings
below show the same cell cut along its two diagonals. In each drawing
the shaded triangle and its unshaded partner give the indicated node.
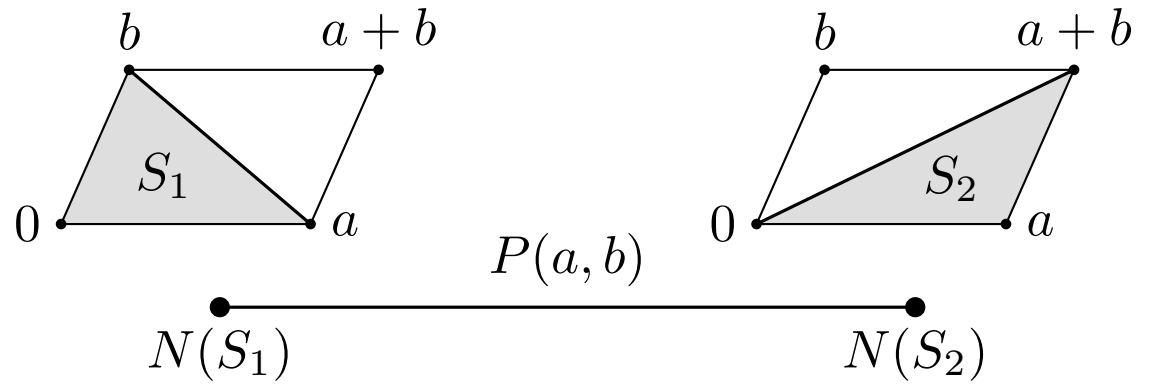
\begin{center}
\begin{tikzpicture}[x=1cm,y=1cm,line join=round,font=\small]
 \foreach \shift/\face in {0/1,4.6/2}{
  \begin{scope}[xshift=\shift cm]
   \coordinate (o) at (-1.05,0);
   \coordinate (a) at (.60,0);
   \coordinate (b) at (-.60,1.02);
   \coordinate (s) at (1.05,1.02);
   \ifnum\face=1
    \fill[black!13] (o)--(a)--(b)--cycle;
    \draw[semithick] (a)--(b);
    \node at (-.37,.32) {$S_1$};
   \else
    \fill[black!13] (o)--(a)--(s)--cycle;
    \draw[semithick] (o)--(s);
    \node at (.24,.30) {$S_2$};
   \fi
   \draw (o)--(a)--(s)--(b)--cycle;
   \foreach \p/\where/\lab in {o/left/0,a/right/a,b/above/b,s/above/{a+b}}
      \fill (\p) circle (1pt) node[\where] {$\lab$};
  \end{scope}}
 \draw[semithick] (0,-.55)--node[above] {$P(a,b)$} (4.6,-.55);
 \fill (0,-.55) circle (1.9pt) node[below] {$N(S_1)$};
 \fill (4.6,-.55) circle (1.9pt) node[below] {$N(S_2)$};
\end{tikzpicture}
\end{center}
The edges for $P(a,b)$ and $P(a,a+b)$ thus meet at $N(S_2)$. In
general, two flat cell orbits are neighbours exactly when their edges
meet at a node. The two nodes of $P(a,b)$ have side vectors $a,b,b-a$ and $a,b,a+b$, respectively,
up to sign. They are distinct, and their two common side directions
determine the parallelogram up to translation, since its sides are
primitive lattice vectors.

\emph{Why each component is a path.}
First, at $N(S)$ every incident flat cell contributes one face of
class $[S]$. This triangle class has just two incidences in the
quotient complex by Lemma~\ref{lem:topology}(2). Hence at most two
actual flat cells meet at any node, so a component cannot branch.

To rule out a circuit, give each node a weight, the sum of the squared
lengths of the three sides of its triangle. We compare this weight
at the two ends of an edge. To pass from one
node to a neighbour, complete the triangle to a parallelogram at one
of its three corners. If $u,v$ are the side vectors pointing out of
that corner, the other node keeps the sides $u,v$ and replaces the
diagonal $v-u$ by $u+v$. The change in the sum is
\[
                  \|u+v\|^2-\|v-u\|^2=4\langle u,v\rangle.
\]
A neighbour therefore has weight no greater than the given node only
when the chosen corner is right or obtuse. A triangle has at most one
such corner, and that corner determines the parallelogram up to
translation. On a circuit, a node of largest weight would have two
such neighbours. This rules out circuits.
There are finitely many cell orbits, so each component of actual flat
cells is a finite path. Its length is the number of edges, one for
each flat cell orbit.

\emph{The end triangles.}
At an interior node of the path, the two flat cells supply opposite
incidences at each of the two triangle classes. These faces cancel in
the fundamental $3$-cycle. At each endpoint of the path,
only one flat cell supplies the two triangles. Each still needs a
second incidence, which must come from a solid neighbour since a flat
neighbour would continue the path. Such a solid cell is called a
\reviewterm{cap}. The two ends give four distinct triangle
classes and four distinct caps up to translation, since a solid
tetrahedron has at most one face parallel to the path plane.

\subsection{Each end cap has index one}
Assume $\theta<2$ and fix a cap at an end triangle. Put $D=\covol(L)$. In these coordinates a lattice vector is short if
it belongs to $T-T$. Write $\gamma$ for the gauge of $T-T$, so $v\in s(T-T)$ exactly
when $\gamma(v)\le s$ for $s\ge0$. By homogeneity, \eqref{eq:rootdet}
gives
\begin{equation}\label{eq:flat-short-triple}
 |\det(x,y,z)|\le 6\vol(T)\gamma(x)\gamma(y)\gamma(z).
\end{equation}
Let $r(P)$ be the least $s\ge0$ such that $P\subset x-sT$ for some $x$.
For a single edge, $0,v\in x-sT$ means that $0=x-sp$ and $v=x-sq$
for some $p,q\in T$. Such an $x$ exists precisely when
$v=s(p-q)\in s(T-T)$. Thus $r(\{0,v\})=\gamma(v)$.

Now consider the four vertices of $P(a,b)$. In coordinates taking $T$
to the standard tetrahedron up to translation, formula~\eqref{eq:example-fit}
gives $r(P)$ as the sum of four maxima. For each of the four linear
facet forms $f$ we have
\[
 \max_{v\in P(a,b)} f(v)
 =\max\{0,f(a),f(b),f(a)+f(b)\}
 =\max\{0,f(a)\}+\max\{0,f(b)\}.
\]
Summing over the four facets gives
$r(P(a,b))=r(\{0,a\})+r(\{0,b\})=\gamma(a)+\gamma(b)$.
By Corollary~\ref{cor:scarf-radius} the radius of each contact cell
is at most the covering radius, which is at most one since $T+L$ covers.

Choose the origin and the edge directions of the end parallelogram so
that the cap has vertices $0,a,b,w$. The scale bound for the
parallelogram gives $\gamma(a)+\gamma(b)\le1$, and the edge
$\{0,w\}$ of the cap gives $\gamma(w)\le1$. Hence the cap index satisfies
\[
       1\le k=|\det(a,b,w)|/D
        \le 6\theta \gamma(a)\gamma(b)\gamma(w)\le 3\theta/2<3.
\]
To eliminate $k=2$, we use the fact that each closed container
contains just its four contacts. Number the facet forms $f_0,f_1,f_2,f_3$ so that their
respective contacts on the parallelogram are $0,a,b,a+b$. Write
$a_i=f_i(a)$ and similarly for other vectors. The signs of
$(a_i,b_i)$ are $(-,-),(+,-),(-,+),(+,+)$ in order. Interchanging
$a,b$ and their facet names if necessary gives $a_3>b_3$ by genericity.
On the three vertices of the shared triangle, the maxima of these
four functions are attained at $0,a,b,a$, respectively.
Facet~3 loses $a+b$ and takes $a$ on the triangle. By
Lemma~\ref{lem:transfer} facet~1 must then pass from $a$ to $w$,
so the cap contacts are $0,w,b,a$. The contact on facet $i$ gives
the upper bound on $f_i(x)$ in its container. The two rows list these bounds.
\[
\begin{array}{c|cccc}
 & f_0(x) & f_1(x) & f_2(x) & f_3(x) \\ \hline
 \text{flat cell's container} & 0 & a_1 & b_2 & a_3+b_3 \\
 \text{cap's container}       & 0 & w_1 & b_2 & a_3
\end{array}
\]
Subtract $f_i(x)$ from the corresponding bound. The difference is
positive on the inner side of the supporting plane, zero on the plane,
and negative on the outer side. Four positive differences put $x$ in
the interior.

Since $a,b$ generate $L\cap\operatorname{span}_{\R}\{a,b\}$, $k=2$ would put
$w$ two lattice layers from that plane. Modulo $2L$, the vector
$w$ is therefore congruent to one of $0,a,b,a+b$. In the first three
cases a cap edge has a lattice midpoint. Hence $w=a+b+2h$ for
$h\in L$ in the intermediate layer. The three points $h+a+b,h+a,h+b$
in this layer must lie outside the closed containers. Their facet
inequalities will put $-h$ inside the flat container.

The cap comparisons $w_1>a_1$ and $w_3<a_3$ imply $h_1>0$ and $h_3<0$.
The point $h+a+b$ is the midpoint of $w$ and $a+b$. It lies strictly
below the bounds on facets $0,2,3$ of the flat cell and $0,1,2$
of the cap. To escape both closed containers it must satisfy
\[
                        h_1+b_1>0,\qquad h_3+b_3>0 .
\]
Consider $h+a$ in the four inequalities for the cap. On facet~1 the
upper bound is $w_1=f_1(w)$. Since $w=a+b+2h$,
\[
 w_1-f_1(h+a)=w_1-(h_1+a_1)=h_1+b_1>0.
\]
Thus $h+a$ strictly satisfies the inequality for facet~1. On facet~2
we have $2(h_2+a_2)=w_2+a_2-b_2<a_2<0$, so $f_2(h+a)<0<b_2$.
On facet~3 we have $f_3(h+a)=h_3+a_3<a_3$. To be outside the
closed cap container, $h+a$ must therefore violate the inequality
for facet~0, giving $h_0+a_0>0$.

Next $h+b=w-(h+a)$ gives $f_0(h+b)=w_0-(h_0+a_0)<0$.
On facets~1 and~3 we have $h_1+b_1=w_1-(h_1+a_1)<w_1$ and
$h_3+b_3<b_3<a_3$. Hence $h+b$ must exceed the bound $b_2$ on
facet~2, giving $h_2+b_2>b_2$ and thus $h_2>0$.

Finally test $-h$ against the four bounds of the flat cell's container.
From $h_0+a_0>0$ we have $-h_0<a_0<0$. Since $h_1,h_2>0$,
we have $-h_1<0<a_1$ and $-h_2<0<b_2$. On facet~3,
$h_3+b_3>0$ gives $-h_3<b_3<a_3+b_3$. Thus $-h$ satisfies all
four inequalities strictly and is a lattice point in the interior,
a contradiction. Hence $k=1$ and $(a,b,w)$ is a basis of $L$.

\subsection{Two bounds on a flat path}\label{sec:analytic-flat-run}
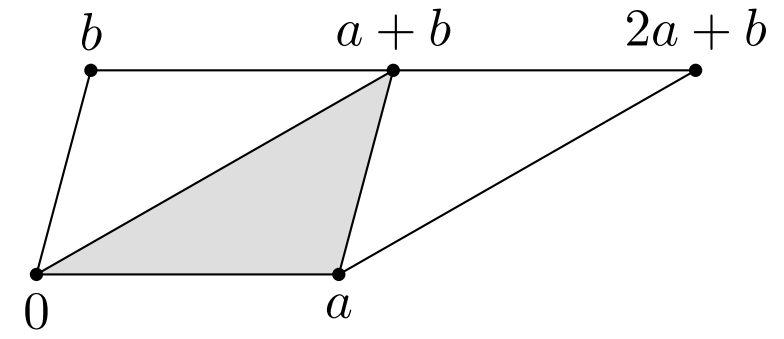
\begin{figure}[htb]
\centering
\begin{tikzpicture}[x=2cm,y=1.35cm,every node/.style={font=\small}]
 \coordinate (o) at (0,0);
 \coordinate (a) at (1,0);
 \coordinate (b) at (.18,1);
 \coordinate (s) at (1.18,1);
 \coordinate (t) at (2.18,1);
 \fill[black!13] (o)--(a)--(s)--cycle;
 \draw (o)--(a)--(s)--(b)--cycle;
 \draw (o)--(s)--(t)--(a);
 \foreach \p in {o,a,b,s,t} \fill (\p) circle (1.25pt);
 \node[below] at (o) {$0$};
 \node[below] at (a) {$a$};
 \node[above] at (b) {$b$};
 \node[above] at (s) {$a+b$};
 \node[above] at (t) {$2a+b$};
\end{tikzpicture}
\caption{The move from $P(a,b)$ to $P(a,a+b)$. The shaded triangle is a face of both cells.}
\label{fig:flat-move}
\end{figure}

Orient a longest path from one of its ends and choose a cap there.
Since $\theta<2$, the three edges of this cap at zero give a basis
$a,b,c=w$ of $L$. The first parallelogram joins the diagonal $b-a$ to $a+b$. At the next node
a continuation must preserve one of the sides $a,b$ and use $a+b$
as its other side. Each subsequent step likewise has one of two forms,
\[
           \mathsf{L}(u,v)=(u,u+v),\qquad \mathsf{R}(u,v)=(u+v,v).
\]
Reading the successive moves along the path gives a \reviewterm{word}
in $\mathsf L$ and $\mathsf R$. A word with $p$ letters describes $p+1$
parallelograms, including the initial one. The empty word describes a single flat cell. Figure~\ref{fig:flat-move} shows the first move $\mathsf{L}$.
Interchanging $a$ and $b$ lets us start every nonempty word with $\mathsf{L}$.
A \reviewterm{run} is a maximal consecutive block of the same letter.

Let $sa+tb$ be the final diagonal, where $s,t$ are positive integers.
The two end diagonals and the cap edge $c$ are short, and
$|\det(a,b,c)|=D$. Thus~\eqref{eq:flat-short-triple} gives
\[
               s+t=|\det(b-a,sa+tb,c)|/D\le6\theta<12.
\]
Hence $s+t\le11$. This bound still permits six consecutive
$\mathsf{L}$ moves from $(a,b)$, whose final diagonal is $7a+b$.
Those six moves give seven cells, with edge pairs $(a,b+ja)$ for
$0\le j\le6$. The following density estimate applies to any generic
covering containing these cells. Its proof uses the common edge $a$
to bound the volume of a set that represents every lattice coset.

\begin{proposition}\label{prop:analytic-flat-run}
Suppose a generic tetrahedron lattice covering contains seven consecutive
flat cells represented, up to translation, by $P(a,b+ja)$ for
$0\le j\le6$. Then its density is at least $256/127$.
\end{proposition}

\begin{proof}
Use barycentric coordinates in which the covering tetrahedron is
$T=\{x\in\R^4: x_i\ge0,\ \sum_i x_i=1\}$. Translation vectors
form the lattice $L$ in the subspace of coordinate sum zero.
Write $s_+=\max\{s,0\}$. Then $\gamma(v)=\sum_i(v_i)_+=\sum_i(-v_i)_+$.
Measure volume in units of $\vol T$ and set $b_j=b+ja$,
$A=\gamma(a)$ and $B_j=\gamma(b_j)$.
For each contact cell $P_j=P(a,b_j)$, the scale formula of the preceding
subsection gives $r(P_j)=A+B_j$. By Corollary~\ref{cor:scarf-radius}
this is at most the covering radius, which is at most one. In coordinate
$i$ the signs of $a_i$ and $(b_j)_i$ determine which vertex attains the
maximum. Each of the four vertices is the unique contact on its facet,
so the four sign pairs $(\operatorname{sgn}a_i,\operatorname{sgn}(b_j)_i)$
all occur.

\emph{The common edge.}
On the two rows with $a_i>0$, both $(b_j)_i$ increase with $j$ and
exactly one is positive for every $j$. Its index is therefore constant.
The same argument applies on the two negative rows. On a row with equal
signs this gives $6|a_i|\le |(b_6)_i|\le B_6$, and on a row with
opposite signs it gives $6|a_i|\le|b_i|\le B_0$. Hence
\[
 |a_i|\le U=\frac{1-A}{6},\qquad A\le 2U,\qquad A\le\frac14 .
\]

\emph{Representatives of the lattice cosets.}
Put $F_+=T\setminus(T+a)$ and $F_-=T\setminus(T-a)$. Each set contains
a representative of every lattice coset, since the representatives of a
coset along any orbit in direction $a$ that meets $T$ have two end
points. For $F=F_+$ or $F_-$, convexity of $T$ gives
$F\cap(F+ka)=\varnothing$ for every nonzero integer $k$.
Since $b_j-b_k=(j-k)a$, the sets $F\cap(F+b_j)$ are therefore
pairwise disjoint. Removing all
seven from $F$ preserves a representative of every coset. Indeed, take a
linear functional $\ell$ with $\ell(a)=0$, $\ell(b)=1$ and choose a
representative in $F$ of least $\ell$-value. Such a representative exists
since $F$ is bounded, and it belongs to none of the $F+b_j$.

\emph{The remaining volume.}
The relative volume removed in passing from $T$ to either $F$ equals
$(1-A)^3$.
To compute the second removal, recall that for a finite set $P$ of
translation vectors,
\[
 \frac{\vol\bigl(\bigcap_{v\in P}(T+v)\bigr)}{\vol T}
       =\left(1-\sum_i\max_{v\in P}v_i\right)_+^3.
\]
Fix $j$, put $\rho=1-B_j\ge A$ and
$\delta_i=\min(|a_i|,|(b_j)_i|)$. Deleting the corner maximizing row
$i$ from $P_j$ reduces its sum of maxima by $\delta_i$. In
$F_+\cap(F_++b_j)$ inclusion--exclusion uses the two triples of
corners containing $0,b_j$. For $F_-$, translate the intersection by
$a$ to obtain the two triples containing $a,a+b_j$. Thus the average
of the volumes of these two intersections is
\[
 \rho^3-\frac12\sum_{i=1}^4(\rho-A+\delta_i)^3+(\rho-A)^3 .
\]

If $A\le 1/5$, the first removal already has volume at least
$64/125>129/256$. Otherwise $1/5<A\le1/4$ and $A/2\le U\le A$.
The two positive values of $a_i$ sum to $A$, as do the absolute values
of the two negative ones. Each is at most $U$. As $t\mapsto(\rho-A+t)^3$
is increasing and convex for $t\ge0$, the two cubes in each pair sum
to at most $(\rho-A+U)^3+(\rho-U)^3$. The average removed for this $j$ is
at least
\[
 \rho^3-(\rho-A+U)^3-(\rho-U)^3+(\rho-A)^3
       =3U(A-U)(2\rho-A)\ge3AU(A-U).
\]
Including the first removal, the total volume removed on average is at
least
\[
 f(A)=(1-A)^3+\tfrac7{12}A(1-A)(7A-1).
\]
On $[1/5,1/4]$ one has $12f'(A)=-43+184A-183A^2<0$, whence
$f(A)\ge f(1/4)=129/256$. At least one of the two remaining sets has
relative volume at most $127/256$. Its lattice translates cover, so
$\covol(L)/\vol T\le127/256$, or equivalently $\theta\ge256/127$.
\end{proof}

\subsection{The possible longest paths}
Return to coverings of density $\theta<2$. By
Proposition~\ref{prop:analytic-flat-run}, six consecutive $\mathsf{L}$
moves anywhere along a path would give $\theta\ge256/127>2$.
Interchanging the two edges gives the same bound for $\mathsf{R}$.
Thus each run has length at most five. We now combine this restriction
with $s+t\le11$.

For a vector $xa+yb$ in the path plane write $\sigma(xa+yb)=x+y$.
If $(u,v)$ is the final edge pair, its diagonal is $u+v=sa+tb$,
so $s+t=\sigma(u)+\sigma(v)$. With one run $\mathsf{L}^p$ the
edge pair is $(a,b+pa)$, giving $s+t=p+2$. All five choices
$1\le p\le5$ meet the diagonal bound.

\emph{Two runs.}
For $\mathsf{L}^p\mathsf{R}^j$, put $v=b+pa$ after the first run.
Each $\mathsf{R}$ move then adds $v$ to the first edge and fixes $v$.
The final pair is $(u,v)$ with $u=a+jv$. Since $\sigma(v)=p+1$ and
$\sigma(u)=1+j(p+1)$, we obtain
\[
       s+t=1+j(p+1)+(p+1)=1+(p+1)(j+1).
\]
Thus $(p+1)(j+1)\le10$. Since $j\ge1$, this forces $p\le4$.
For $p=1,2,3,4$ the largest possible $j$ is respectively $4,2,1,1$,
giving eight words.

\emph{Three runs.}
Append $\mathsf{L}^r$ to the same two runs. The first edge $u$ now
stays fixed while the second becomes $v+ru$. Their sum is $v+(r+1)u$,
so the bound becomes
\[
       s+t=p+1+(r+1)\bigl(1+j(p+1)\bigr)\le11.
\]
Since $p,r\ge1$, the left side is at least $2+2(1+2j)$, forcing
$j=1$. The bound is then $p+1+(r+1)(p+2)\le11$. For $p=1$
it gives $r\le2$, and for $p=2$ it gives $r=1$. Larger $p$ are
excluded already at $r=1$. The three pairs are $(1,1),(1,2),(2,1)$.

\emph{Four or more runs.}
The pair of coefficient sums $(\sigma(u),\sigma(v))$ follows the same
addition rule as $(u,v)$. Each move increases one entry and both moves
preserve coordinatewise inequalities. Thus shortening each of the first
four runs to one move and omitting all later moves gives a lower bound
for the final sum. For $\mathsf{LRLR}$ the successive pairs are
\[
 (1,1)\xrightarrow{\mathsf L}(1,2)\xrightarrow{\mathsf R}(3,2)
 \xrightarrow{\mathsf L}(3,5)\xrightarrow{\mathsf R}(8,5).
\]
Hence four runs already give $s+t\ge8+5=13$, exceeding eleven.
Including the empty word, we have $1+5+8+3=17$ words before reversal.

Reading a path from its other end reverses its word. Indeed, write
the last parallelogram $P(u,v)$ as the translate $P(-u,v)+u$.
In these coordinates the reverse of a move has the same letter, since
\[
 \mathsf{L}(-u,u+v)=(-u,v),\qquad
 \mathsf{R}(-(u+v),v)=(-u,v).
\]
Swapping the two initial edges interchanges $\mathsf{L}$ and $\mathsf{R}$.
Reversing a word and, when necessary, swapping the edges gives
\[
 \mathsf{LLR}\sim \mathsf{LRR},\qquad \mathsf{LLLR}\sim \mathsf{LRRR},\qquad
 \mathsf{LLRL}\sim \mathsf{LRLL},\qquad \mathsf{LLLLR}\sim \mathsf{LRRRR}.
\]
Up to these symmetries, the sequence of flat cells along a longest
path is described by one of the thirteen words
\begin{equation}\label{eq:flat-words}
 \begin{gathered}
 \varnothing,\quad \mathsf{L},\quad \mathsf{LL},\ \mathsf{LR},\quad \mathsf{LLL},\ \mathsf{LLR},\ \mathsf{LRL},\\
 \mathsf{LLLL},\ \mathsf{LLLR},\ \mathsf{LLRL},\ \mathsf{LLRR},\quad \mathsf{LLLLL},\ \mathsf{LLLLR} .
 \end{gathered}
\end{equation}
The ancillary command \texttt{python3 examples.py words} reproduces the
list by applying the two moves and identifying reverse paths.

Choose the initial end and its two plane edges to give the selected
representative. If this reverses the path, take an actual cap at the
new initial end. Its opposite vertex again completes the plane basis
by the height argument. In these coordinates write the basis as
$(a,b,c)$ and number the facets so that the contacts on $P(a,b)$
are $(a+b,a,b,0)$.
After either move the forms positive on both edges or negative on both
keep those signs. The two other forms have opposite signs on the
persisting edge. Four distinct contacts require each to take the
opposite sign on the new edge. Hence every edge pair $(u,v)$ of the
word gives the contacts $(u+v,u,v,0)$ in facet order. Undoing the possible
interchange in the cap argument gives the two facet assignments for the cap,
\[
                       (a,c,b,0),\qquad (b,a,c,0).
\]
For the empty word, swapping $a,b$ and the two middle facets identifies
the two choices on the entire complex. For each other word we retain
both. This proves Proposition~\ref{prop:flat-catalogue}.

For generic coverings of density less than two,
Propositions~\ref{prop:solid-catalogue} and~\ref{prop:flat-catalogue}
answer the question in the introduction about four facet contacts
forming an affine basis of $L$. On passing to a limit, the four points
remain on their respective closed facets.

\begin{corollary}\label{cor:unit}
Suppose $T+L$ is a tetrahedron lattice covering of density less than
two. For some $0<r\le1$ a translate of $-rT$ has lattice-free
interior and contains, on its four respective closed facets, points
$p_1,p_2,p_3,p_4\in L$ forming an affine basis of $L$. For a generic
covering the four points can be chosen as the contacts of a maximal cell.
\end{corollary}
\begin{proof}
For a generic covering take the cell of index one from
Proposition~\ref{prop:solid-catalogue} if all cells are solid,
and a cap of a flat component otherwise. In either case its scale
is at most one by Corollary~\ref{cor:scarf-radius}.

In the general case, normalize to $T=\Delta$ and $L=B\Z^3$.
Lemma~\ref{lem:perturb} gives generic covering bases $B_j\to B$
of density less than two. For each $j$ translate the first contact
of the index-one cell to zero. The integer labels then belong to
$B_j^{-1}(\Delta-\Delta)$, a uniformly bounded family. Along a subsequence the four labels, listed in facet order,
are constant and their edge determinant is $\pm1$ in
lattice coordinates. Their support values converge to those of the same
labels at $B$. The limiting scale is at most one and is positive by
affine independence. Each contact remains on its assigned facet. An
integer label in the interior would satisfy all four strict inequalities
also for large $j$, contrary to lattice emptiness there.
\end{proof}

To continue from a longest path in a generic covering of density
$\theta<2$, use the basis $a,b,c$ of its chosen initial cap.
For a short vector $v=xa+yb+zc$, we have $\gamma(c)\le1$
and $\gamma(a)+\gamma(b)\le1$. Hence~\eqref{eq:flat-short-triple} gives
\[
\begin{aligned}
 D(|x|+|y|)&=|\det(v,b,c)|+|\det(a,v,c)|
     \le 6\vol(T)\gamma(v)\gamma(c)(\gamma(a)+\gamma(b))\le 6\vol(T),\\
 D|z|&=|\det(a,b,v)|\le6\vol(T)\gamma(a)\gamma(b)\gamma(v)\le\tfrac32\vol(T).
\end{aligned}
\]
Thus $|x|+|y|\le11$ and $|z|\le2$. Relative to any vertex of an
unpaired face this includes every possible fourth vertex. Further
short edges impose the same determinant test. A flat neighbour must be
a fundamental parallelogram, and the height argument applies to every
solid neighbour at a flat face. Its facet assignment must satisfy
Lemma~\ref{lem:transfer}. Since the chosen path was longest, a component
longer than it may also be discarded in the continuation. These are the
initial configurations and restrictions used in Section~\ref{sec:finite}.

\section{Determinants on the remaining supports}\label{sec:finite}

We shall separate the support of the $84$-shape from the others at
the covolume $D_0=83988/10^6$, just below $21/250$.
\begin{proposition}\label{prop:terminal}
Every generic covering $\Delta+B\Z^3$ has $D=|\det B|\le21/250$.
If $D>D_0$, then, after a change of lattice basis and a permutation
of the facets, its support consists of the thirteen cells
in~\eqref{eq:84-cycle}, with the integer labels of
Appendix~\ref{sec:sharp-identity}. In these coordinates $\det B<0$.
\end{proposition}
Thus a generic covering on any other support has density at least
$1/(6D_0)=125000/62991>125/63$.

\subsection{Inequalities from a list of contacts}\label{sec:checks}
Suppose that $D>D_0$. Then $\theta=1/(6D)<2$, so the classifications
of Sections~\ref{sec:noflat} and~\ref{sec:flat-components} apply.
Consider the first parallelogram of a flat path and one of its caps.
As in Section~\ref{sec:flat-components}, take the cap's three edge
vectors at the origin as the columns of the lattice basis matrix
$B=[\,a\ b\ c\,]$. The parallelogram has vertices $0,a,b,a+b$,
and the cap has vertices $0,a,b,c$. In the basis $B$ their vertex
sets in $\Z^3$ are therefore
\[
 P=\{0,e_1,e_2,e_1+e_2\},\qquad Q=\{0,e_1,e_2,e_3\}.
\]
Number the facets as in the proof of Proposition~\ref{prop:flat-catalogue}
and take the first of the two facet assignments for the cap. In
facet order, the contact labels for $P$ and $Q$ are, respectively,
\[
 (e_1+e_2,e_1,e_2,0),\qquad (e_1,e_3,e_2,0).
\]
In each tuple the $i$th entry is the vertex at which $g_i$ attains
its maximum on the corresponding set. For example, the second
entries give $\max_{p\in P}g_2(p)=g_2(e_1)$ and
$\max_{q\in Q}g_2(q)=g_2(e_3)$.

For either vertex set $F$, the least $t\ge0$ for which a translate of
$-t\Delta$ contains $BF$ is $r_B(F)=\sum_{i=1}^4\max_{v\in F}g_i(v)$
by~\eqref{eq:example-fit}. Since $\Delta+B\Z^3$ covers, both radii
are at most one by Corollary~\ref{cor:scarf-radius}. Reading off the four maxima and using
$g_i(e_j)=B_{ij}$ for $1\le i\le3$ and $g_4(0)=0$, we obtain
\[
\begin{aligned}
 r_B(P)&=g_1(e_1+e_2)+g_2(e_1)+g_3(e_2)+g_4(0)\\
       &=B_{11}+B_{12}+B_{21}+B_{32}\le1,\\
 r_B(Q)&=g_1(e_1)+g_2(e_3)+g_3(e_2)+g_4(0)\\
       &=B_{11}+B_{23}+B_{32}\le1.
\end{aligned}
\]
The translates of $-\Delta$ containing $BP$ and $BQ$ may differ.
Both inequalities constrain the same matrix $B$.

More generally, let $S$ be a list of cells with specified facet assignments, in the
basis supplied by an index-one cell (the initial cap in the flat case).
For each vertex set $P$ in $S$, let $(w_1,w_2,w_3,w_4)$ be its contact
labels in facet order. This assignment and
Corollary~\ref{cor:scarf-radius} give
\[
 g_i(w_i-p)>0\quad(p\in P\setminus\{w_i\}),\qquad
 L_P(B):=1-\sum_{i=1}^4g_i(w_i)=1-r_B(P)\ge0 .
\]
Here $g_1,g_2,g_3$ are the rows of $B$, and $g_4$ is their negative sum.
The closed container is the intersection of the four half-spaces
$g_i(y)\le g_i(w_i)$. Its only integer points are $w_1,w_2,w_3,w_4$.
Thus each $x\in\Z^3\setminus P$ lies strictly beyond at least one
supporting plane. This gives the \reviewterm{exit condition}
\begin{equation}\label{eq:or}
 g_i(x-w_i)>0 \qquad\text{for at least one }i\in\{1,2,3,4\}.
\end{equation}
The facet index may depend on $x$.
The conditions are invariant under a common translation of $P$ and $x$.

For facet $i$, the \reviewterm{row cone} $C_i$ consists of the row
vectors $g$ satisfying $g(w_i-p)\ge0$ for every prescribed comparison.
It contains no line, since the three edge differences at $w_i$ in a
solid cell span $\R^3$. We represent each extreme ray (one-dimensional
face) by its primitive integer generator. To find these generators,
take both signs of the cross products of pairs of independent defining
normals, keep the feasible directions and divide by the coordinate gcd. A linear
form can be positive on $C_i$ only if it is positive on some ray.
All strict defining comparisons can hold simultaneously precisely
when each is positive on some ray, since adding all the rays then
gives a strict solution. Such a solution lies in the interior.
For a fixed $x$, rule out facet $i$ if $g_i(x-w_i)$ is nonpositive
on every ray of $C_i$. If this form is nonnegative on every ray and
positive on some ray, facet $i$ already ensures~\eqref{eq:or} for $x$.
If no facet is possible, discard the list of cells. If only one facet
is possible, impose $g_i(x-w_i)>0$ for that row.
Recompute the rays and repeat to stability, rejecting a
cone with no strict solution. Appendices~\ref{sec:unit-runs}
and~\ref{sec:flat-runs} specify the finite sets of labels used at
each stage. For determinant bounds on a paired support, translate the
lexicographically least vertex of each cell to zero and start the row cones anew. Use the
integer labels in the boxes $p+[-2,2]^3$, $p\in P$. For each pair of
distinct vertices $p,p'\in P$ with midpoint $m=(p+p')/2$, also take
the eight vertices $\lfloor m\rfloor+\{0,1\}^3$ of the unit cube with
lower corner $\lfloor m\rfloor$, taking the floor in each coordinate.
Omit the vertices of $P$ from these tests.

Two further inequalities will be used on particular supports.
For the first, fix a nonzero label $v\in\Z^3$, put $y=Bv$ and
write $\gamma$ for the gauge of $\Delta-\Delta$. Every lattice coset
has a representative in $\Delta$. Subtract $y$ as long as the
representative stays in $\Delta$. Boundedness gives a representative in
$\Delta\setminus(\Delta+y)$, so this set has volume at least $D$.
Taking the tighter bound on each of the four facets shows that
$\Delta\cap(\Delta+y)$ is a translate of $(1-\gamma(y))\Delta$ for
$\gamma(y)\le1$ and is empty for $\gamma(y)>1$. Hence
\begin{equation}\label{eq:flat-period}
  6D\le1-(1-\gamma(y))_+^3.
\end{equation}
Since $6D_0>1-(19/24)^3$, this gives $24\gamma(Bv)>5$.
On a cone where the signs of all four $g_i(v)$ are fixed,
$\gamma(Bv)=\sum_i(g_i(v))_+$ is linear in $B$.

For the second inequality, let $H=\conv P$ and $H'=\conv P'$
for two prescribed cells, allowing repetition. Both $BH$ and $BH'$ lie
in translates of $-\Delta$. By the facet support formula,
$6V(-\Delta,-\Delta,[0,y])=\gamma(y)$. Monotonicity, translation invariance,
and a change of variables give
\begin{equation}\label{eq:gauge-mixed}
 D\,6V(H,H',[0,v])
   =6V(BH,BH',[0,Bv])\le \gamma(Bv).
\end{equation}
Taking the maximum over the prescribed cells and using $D>1/12$ gives
\[
     m(v)=\max_{P,P'}6V(\conv P,\conv P',[0,v]),\qquad
     12\gamma(Bv)\ge m(v).
\]
Use one sign of every primitive
$v\in[-3,3]^3\cap\Z^3$ and of every prescribed edge whenever the
four row signs are determined. Appendix~\ref{sec:ray-weights} gives
a planar area formula for $m(v)$.

Together with $1$, the contact comparisons, the forms $L_P$, and the
primitive inward facet forms of the resulting cones, these give a
finite family of nonnegative affine forms in the nine entries of $B$.
The columns $Be_j$ are edges of the index-one cell used to choose the basis and therefore
lie in $\Delta-\Delta$, so $1\pm B_{ij}\ge0$. All affine calculations substitute
$g_4=-g_1-g_2-g_3$. In particular, a nonnegative combination with
coefficients $(c_0,c_1,\ldots,c_9)$ satisfying
$c_0+\sum_{j=1}^9|c_j|<0$ is already a contradiction.

\subsection{Two bounds for the determinant}\label{sec:last}
The first bound uses the determinant's multilinearity to treat both
signs at once. The absolute determinant of any three of the four
rows is $D$ because their sum is zero.

\begin{lemma}[Three-row bound]\label{lem:three-row}
Let $i,j,k$ be three distinct row indices. Suppose $b,w>0$ and
$h_i,h_j,h_k\in\R^3$ satisfy
$g_i(h_i)+g_j(h_j)+g_k(h_k)\le w$ on the prescribed domain.
Suppose $r\cdot h_i,s\cdot h_j,u\cdot h_k$ are nonnegative on
the primitive rays $r,s,u$ of $C_i,C_j,C_k$, respectively, and
\begin{equation}\label{eq:three-rays}
 |\det(r,s,u)|w^3
 \le27b\,(r\cdot h_i)(s\cdot h_j)(u\cdot h_k)
\end{equation}
for every triple of such rays. Then $D\le b$.
\end{lemma}
\begin{proof}
Expand each row as a nonnegative combination of its extreme rays.
Multilinearity and the triangle inequality extend
\eqref{eq:three-rays} to the three rows themselves. The three
factors on the right have sum at most $w$, so their product is at
most $(w/3)^3$. Cancelling $w^3$ gives $D\le b$.
\end{proof}

A nonnegative combination of the affine inequalities supplies the
linear bound required by the lemma. Its weights and the ray comparisons
can be tested with integer arithmetic as described in
Appendix~\ref{sec:ray-weights}. If a cone is cut by a rational plane
through the origin, proving the bound on both closed pieces proves
it on the whole cone.

The second bound treats each determinant sign $e\in\{-1,1\}$
separately. Let $\mathcal L$ be any of the finite lists of
nonnegative forms specified below. A triple product of its members
and any product $LH^2$, with $L\in\mathcal L$ and $H$ affine, are
nonnegative. For a nonnegative rational combination $P$ of such
products write $e\det B+P(B)=\sum_\alpha c_\alpha B^\alpha$.
Since $|B_{ij}|\le1$, each monomial $B^\alpha$ has absolute value at most one. Thus
\begin{equation}\label{eq:cubic-terminal}
 e\det B\le e\det B+P(B)\le c_0+\sum_{\alpha\ne0}|c_\alpha|.
\end{equation}
For a specified finite list of $H$, testing a rational upper bound
$b$ in \eqref{eq:cubic-terminal} is a rational linear feasibility
problem in the product weights and auxiliary variables
$t_\alpha\ge\pm c_\alpha$ with $c_0+\sum_{\alpha\ne0}t_\alpha\le b$.
Expansion after clearing denominators verifies a proposed bound.

\subsection{The terminal cases}
\begin{proof}[Proof of Proposition~\ref{prop:terminal}]
It suffices to consider $D>D_0$, which gives $\theta<2$. If all cells are solid,
Proposition~\ref{prop:solid-catalogue} gives three candidates for the
entire support, with six, six and fourteen cells.
For each support take the constant
$1$, the forms $L_P$, and the facet forms of the tightened row cones.
Appendix~\ref{sec:unit-runs} gives the six signed bounds from
\eqref{eq:cubic-terminal}. Their common consequence is
\[
                   D\le 83396/10^6<D_0.
\]

If a flat cell occurs, choose a longest path.
Proposition~\ref{prop:flat-catalogue} gives thirteen possible words
for this path. Each word specifies a list of parallelograms to which
we adjoin the first cap with one of its allowed facet assignments. For the three
words $\mathsf{LLRL},\mathsf{LLRR},\mathsf{LLLLR}$ these cells alone
give the required bound. For each edge pair $(u,v)$ in the word the facet assignment
gives $\gamma(Bu)=(g_1+g_2)(u)$ and $\gamma(Bv)=(g_1+g_3)(v)$.
Use all triple products of $1$, the forms $L_P$, the contact comparisons,
and $24\gamma(Bu)-5,24\gamma(Bv)-5$ for these pairs. Formula
\eqref{eq:cubic-terminal} gives the upper bounds
\[
  \begin{array}{c|rrr}
     \text{word}&\mathsf{LLRL}&\mathsf{LLRR}&\mathsf{LLLLR}\\ \hline
     10^6D\le &76844&82287&77046
  \end{array}
\]
The command \texttt{python3 bounds/run.py paths} checks the two determinant
signs for both facet assignments on each of the three paths. All give $D<D_0$.

For the remaining ten words, complete the prescribed path and cap
by finding partners for their unpaired faces as in Appendix~\ref{sec:flat-runs}.
The seven words $\mathsf{LR}$, $\mathsf{LLL}$, $\mathsf{LLR}$, $\mathsf{LRL}$,
$\mathsf{LLLL}$, $\mathsf{LLLR}$ and $\mathsf{LLLLL}$ give $3021$
paired lists, which reduce to $945$ supports after a common
change of lattice basis and permutation of the rows. On $944$
supports Lemma~\ref{lem:three-row} applies directly with $b=D_0$. For
representative $65$ it applies on both sides of $g_3(0,2,3)=0$ in the normalized
coordinates. Thus $D\le D_0$ on all $945$.

The empty word, $\mathsf{L}$, and $\mathsf{LL}$ give $1501$ paired lists and $399$
distinct normalized supports. Two have inconsistent affine inequalities.
The three-row bound gives $D\le D_0$ on $347$ directly and on $40$
after subdividing row cones. For nine of the last ten,
\eqref{eq:cubic-terminal} gives $D\le83562/10^6<D_0$.
On representative $104$, the positive sign gives
$\det B\le17708/10^6<D_0$, and the identity in
Appendix~\ref{sec:sharp-identity} gives $-\det B\le21/250$.
Thus only this last support with negative determinant can occur
above $D_0$.
\end{proof}

The thirteen cells of representative $104$ are the terms of the
$84$-shape cycle~\eqref{eq:84-cycle}. Appendix~\ref{sec:sharp-identity}
gives their integer labels in a basis $A$ of $C\Z^3$ and the
coefficient identity for $21/250+\det B$ on this support.
Interchanging the second and fourth facet forms amounts to applying
$U(x,y,z)=(x,-x-y-z,z)$. At $B_*=UA/10$ all thirteen cells have
radius one and $-\det B_*=21/250$, attaining the cubic bound.

\subsection{The optimal density}\label{sec:sharp}
\begin{proof}[Proof of Theorem~\ref{thm:main}]
The degree-one cycle~\eqref{eq:84-cycle} gives the covering
$\Delta+L_{84}$ of density $125/63$. Its lattice has covolume
$84/10^3=21/250$.

Let $\Delta+B\Z^3$ be any covering. By Lemma~\ref{lem:perturb} it is
the limit of generic coverings with bases $B_j\to B$.
Proposition~\ref{prop:terminal} gives $|\det B_j|\le21/250$ and hence
$|\det B|\le21/250$.

Suppose $|\det B|=21/250$. For all large $j$, $|\det B_j|>D_0$,
so the same proposition puts each covering on the thirteen-cell support
of the $84$-shape. In its normalized coordinates the basis matrix is
$S_jB_jV_j$, where $V_j\in\GL_3(\Z)$ selects the three edges at a
vertex of an index-one cell. The map $S_j$ permutes the four facet
forms. These maps are the linear parts of the affine symmetries of
$\Delta$, as one sees by permuting its four barycentric coordinates.
The columns of $B_jV_j$ lie in $\Delta-\Delta$. Since $B_j\to B$ and
$B$ is invertible, the integer columns of $V_j$ lie in a bounded set.
After passing to a subsequence we thus have $V_j=V$ and $S_j=S$.

The limit $SBV$ has determinant $-21/250$. Each summand on the right
of the identity in Appendix~\ref{sec:sharp-identity} is nonnegative on
the normalized matrices and remains so in the limit. The equality
case of that identity gives $SBV=B_*$. Since $B_*\Z^3=U L_{84}$,
the lattice $B\Z^3=S^{-1}U L_{84}$ is obtained by a symmetry of
$\Delta$. Conversely, the linear part of any such symmetry carries
$L_{84}$ to another optimal covering lattice for $\Delta$. An affine
map taking $\Delta$ onto $T$ transports its covering lattices by its
linear part and preserves their densities. This proves the assertion
for $T$.
\end{proof}

\begin{proof}[Proof of Corollary~\ref{cor:cayley}]
Send $e_1,e_2,e_3$ to the three generators and let $L$ be the kernel
of the resulting homomorphism from $\Z^3$ onto the group. A directed
path of length at most $d$ gives a nonnegative representative of
coordinate sum at most $d$ in each coset. For $x=p+u$ with $p\in\Z^3$
and $u\in[0,1)^3$, adding $u$ to such a representative of $p+L$
shows that $x\in L+(d+3)\Delta$. Theorem~\ref{thm:main} now gives
\[
               [\Z^3:L]=\covol(L)\le \frac{21(d+3)^3}{250}.
\]

For a positive integer $k$, put $L=kC\Z^3$. Scaling the covering
given by~\eqref{eq:84-cycle} yields $10k\Delta+L=\R^3$.
Given $p\in\Z^3$, apply this covering to $p+u$ with
$u=(3/4,3/4,3/4)$. Choose $\lambda\in L$ such that
$p+u-\lambda\in10k\Delta$. Each coordinate of $p-\lambda$ is an
integer at least $-3/4$, hence is nonnegative. Their sum is an
integer at most $10k-9/4$, hence at most $10k-3$. Thus every coset
has a nonnegative representative of coordinate sum at most $10k-3$,
and the directed diameter is at most $10k-3$.
The group has order $84k^3$. Substituting this
order in the bound just proved gives $d\ge10k-3$, so equality holds.
\end{proof}

\appendix
\section{Finite checks for the cell indices}\label{sec:minindex}

\subsection{Excluding high-index cells}\label{sec:high-replay}
Suppose a generic covering of density $\theta<2$ contains a cell with
hull $H=H_{m,t}$, where $6\le m\le9$. In the contact cycle every
triangle belongs to two cells. Start with $H$, choose one of its
triangular faces, and ask which cell could meet $H$ there. Each possible
neighbour, together with its facet assignment, gives a separate choice
to pursue. After each addition we ask whether the cells chosen so far
can still belong to a closed cycle of odd degree.

\emph{The candidate cells.}
Use the lattice coordinates in which $H$ has the given form. Translating
the lexicographically least vertex of a cell to zero gives one vertex
set per translation class. By~\eqref{eq:highwindow} its four labels
satisfy $a(v)\in[-11,11]^3$. Among the four-point sets in this window
satisfying~\eqref{eq:highwindow}, retain all solid tetrahedra of index
at most nine whose only lattice points are their vertices. In the search
allowing flat cells at $m=9$, also retain the coplanar sets. Every cell
of the covering has a representative on the resulting list.

\begin{samepage}
\emph{One step across a triangle.}
For the all-solid search with $H=H_{7,1}$, put $q=(1,7,1)$. Number the
facets so that their contacts on $H$ are in lexicographic vertex order.
Take the face $F=\{e_3,e_1,q\}$ opposite $0$ in $H$. A translate of one candidate $P$ has the same
three labels on a face. For each facet of the two containers, the table
gives the vertex that touches it.
\[
\begin{array}{c|cccc}
 &\text{facet }1&\text{facet }2&\text{facet }3&\text{facet }4\\ \hline
 H&0&e_3&e_1&q\\
 P&e_3&-e_2&e_1&q
\end{array}
\]
\end{samepage}
The two cells share the three vertices of $F$. Across $F$, the shared
vertex $e_3$ takes over facet~1 from $0$, and the new vertex $-e_2$
takes over facet~2 from $e_3$. The other two contacts keep their facet
numbers, as required by Lemma~\ref{lem:transfer}. With the facet
assignments in the table, the contact and mixed-volume tests described
below leave $265$ of the $432$ initial candidate vertex sets. If a
covering realizes this choice, its contact cycle must solve the face and
degree equations~\eqref{eq:cup} on this reduced list with $x_H=x_P=1$.
This system has no solution, ruling out the proposed neighbour with
these contacts.

There are $27$ choices at $F$, counting each neighbour with its
assignment. For each choice we restrict~\eqref{eq:cup} to the vertex
sets surviving its contact and mixed-volume tests and fix the variables
of the selected cells at one. The resulting system
has no solution for $22$ of the $27$ two-cell selections, including the
one displayed above. The other five give $22$ three-cell selections,
all of which the same test rules out.
Each choice is a node of the search tree. Its children are the choices
obtained by prescribing one more cell and its facet assignment.
This accounts for the $1+27+22=50$ nodes visited for $(7,1)$ and
exhausts every possible continuation from $H$.

\emph{Testing and extending a choice.}
Write $S$ for the cells selected so far, each with its assigned contacts.
In a covering the containers of all cells use the same four facet
functionals. As above, number the facets by their contacts on $H$.
If $h$ is its contact on facet $i$, the possible functionals $g$ for
this facet satisfy $g(h)>g(w)$ at every other vertex $w$ of $H$.
These directions form the interior of the normal cone at $h$.
For two vertices $v,w$ of a candidate, the order of $g(v)$ and $g(w)$
can change only across the plane $g(v-w)=0$. The planes for all such
pairs divide the cone into \reviewterm{chambers}. Within a chamber
the maximizing vertex on each candidate is fixed. Genericity puts the
actual $g_i$ in one of these chambers. For each facet $i$ retain just
those chambers giving the prescribed contact on every cell in $S$.
A new assignment must use four distinct vertices, and for each $i$
some retained chamber must give its assigned vertex as maximizer.

Also remove candidates that violate the mixed-volume bound of
Lemma~\ref{lem:mix} with the selected cells. To test whether $S$ could
still be completed, set $x_P=1$ for $P\in S$ in~\eqref{eq:cup}.
Each other surviving candidate may be used ($x_P=1$) or omitted
($x_P=0$). If the face and degree equations have no common solution,
every completion of this choice is ruled out.

For a selection passing these tests, count its triangular faces up to
translation. A face occurring once still needs a neighbour. At each
such face list every unselected cell and contact assignment which,
after aligning the face, give the transfer of Lemma~\ref{lem:transfer}
illustrated above. In facet order this gives opposite orientations
on the common face. The search proceeds as follows.
\begin{quote}\begin{minipage}{\linewidth}\small\ttfamily\raggedright
Start with $S$ consisting of $H$ and its assigned contacts.\\
At each node, update the contact chambers and apply the mixed-volume bound.\\
Reject if a selected cell is removed or \eqref{eq:cup} has no solution.\\
Reject a third incidence at a face, an incompatible pair, or an unpaired face with no possible neighbour.\\
If every face is paired or nine cells have been added to $H$, count this node as unresolved and return.\\
Choose an unpaired face with the fewest choices. For each choice on its list, add that cell and assignment to $S$ in a separate branch and repeat.
\end{minipage}\end{quote}
If a covering realizes $S$, its facet functionals realize the assigned
contacts, its cells satisfy Lemma~\ref{lem:mix}, and its full cycle
solves~\eqref{eq:cup}. At whichever unpaired face we choose,
its actual neighbour belongs to a different translation class and
satisfies the transfer rule (Lemmas~\ref{lem:topology}(2)
and~\ref{lem:transfer}). One of the branches therefore continues to
follow that covering. It must eventually reach a node with all faces
paired or with nine cells added to $H$. Either outcome counts as
unresolved. Thus zero unresolved nodes exclude a covering containing $H$.

The command \texttt{python3 indices/run.py high} runs the all-solid
search for each of the seven pairs $(m,t)$ in the proof of
Proposition~\ref{prop:high}. When flat cells are allowed the bound to
prove is $m\le8$, so only $(9,1)$ and $(9,2)$ require a further search.
The table gives the initial number of candidate cells up to translation
and the number of nodes visited. For $(8,1)$ there
are $210$ solid candidates, but the face and degree equations already
have no common solution with $x_H=1$. Only the initial node is visited.
\begin{center}
\begin{tabular}{ccrr|rr}
$m$&$t$&\multicolumn{2}{c|}{solid cells only}&\multicolumn{2}{c}{flat cells allowed}\\
&&candidates&nodes&candidates&nodes\\ \hline
6&1&668&248&&\\
7&1&432&50&&\\
7&2&452&140&&\\
8&1&210&1&&\\
8&3&207&54&&\\
9&1&94&1&162&59\\
9&2&80&1&129&53
\end{tabular}
\end{center}
All nine searches have zero unresolved nodes. This excludes index nine
even with flat cells and indices six through nine in the all-solid case,
as required for Proposition~\ref{prop:high}.

\subsection{Excluding minimum solid indices two, three and four}\label{sec:minimum-indices}
Suppose a generic covering has only solid cells and $\theta<2$. In the coordinates of
Section~\ref{sec:cell-index-one}, a cell of minimum index $q$ is
$V=\{0,qe_1,qe_2,qe_3\}$ in $\Gamma_q$. Translate the least vertex of
each candidate to zero. Its hull must be empty of other lattice points,
its index must lie in $[q,5]$, and its six edge differences must be
primitive in $\Gamma_q$ and belong to $U_q$. In particular the finite
set $U_q$ contains every possible fourth vertex at a triangle translated
to have a vertex at zero. By~\eqref{eq:noflat-triples}, any three edges
from cells of the covering have absolute determinant at most $11q^2$.

\emph{Minimum index three or four.}
Start with $V$ and forget the facet assignments. Every actual covering
still gives a selection of vertex sets. At an unpaired face try all candidates in the window above.
Discard a choice if it gives a third incidence at a face or violates
the edge determinant bound. If a list is empty, reject the branch.
If an unpaired face has only one possible neighbour, that neighbour is
forced. Add all such neighbours. If unpaired faces remain and none has a
unique partner, choose one with the fewest possible neighbours and try each in turn. Each addition introduces a
new cell orbit from a finite list. The actual neighbour across each
unpaired face satisfies all these tests, so a covering yields a branch on
which every face is eventually paired.

A connected set of actual cells with every face paired is the entire
quotient complex, since the dual graph is connected. Orienting its cells
by opposite incidence at shared faces must therefore give signed index
sum $\pm6$, by Lemma~\ref{lem:topology}(3) and~\eqref{eq:cycle-degree}. The command
\texttt{python3 indices/run.py three-four} visits $520$ distinct sets
for $q=3$. Exactly four have every face paired. Each admits a consistent
orientation and has signed index sum zero. For $q=4$
none of the $17$ sets has every face paired. Both minimum indices are
excluded.

\emph{Minimum index two.}
Fix $F=\{0,2e_1,2e_2\}$ in $V$. Enumerating the fourth vertex in
$U_2\setminus V$ subject to the edge, index and emptiness conditions above gives $120$
possibilities for the neighbour across $F$. For each, try the three
transfers of Lemma~\ref{lem:transfer}, one for each vertex of $F$.
These $360$ first pairs with assigned contacts include the neighbour
of $V$ in any actual covering. The command \texttt{python3 indices/run.py first}
generates these assignments.

For each candidate tetrahedron try all assignments of its four vertices
to the four facets, respecting the prescribed contacts. For each row, its prescribed
contacts must strictly maximize that row on their cells. Moreover, for
a prescribed or trial cell $P$ with contacts $w_i$, every
$x\in U_2\setminus P$ must satisfy $g_i(x-w_i)>0$ for some $i$ by
Lemma~\ref{lem:topology}(1). Test these alternatives on the row cones
as in~\eqref{eq:or}. If none is possible, discard the assignment. If only one remains
on a prescribed cell, its inequality can be imposed
on that row. Thus these restrictions preserve the rows of every covering
extending the prescribed cells.

Two assignments can be paired across a common triangle if the cells are
distinct modulo translation, induce opposite orientations on the triangle,
and satisfy the contact transfer of Lemma~\ref{lem:transfer}. Discard assignments incompatible
with prescribed incidences or left without a partner at some face.
In a covering the two actual assignments at each face survive each
round of deletions together. Reject a branch if a prescribed assignment
is lost. At a surviving node either every face is paired or some
unpaired face remains, in which case try all its partners. The command \texttt{python3 indices/run.py two} visits $360$,
$75$ and $6$ nodes after zero, one and two additions, respectively.
All branches are rejected. Grouping the initial pairs under the six
permutations of $F$ gives the same exclusion from $65$ representatives.
This excludes $q=2$ and finishes the minimum-index step of
Proposition~\ref{prop:solid-catalogue}.

\section{The all-solid calculation at index one}\label{sec:unit-runs}
The enumeration starts with one marked face of an index-one cell.
It produces twenty paired lists. We give the changes of coordinates
that identify their three distinct supports and the determinant bound for each sign.
The command \texttt{python3 solid/run.py} generates the first pairs,
completes all branches, and normalizes the paired lists.

\subsection{An index-one cell and its first neighbour}
Take the vertex set $V=\{0,e_1,e_2,e_3\}$ and the marked face
$F=\{0,e_1,e_2\}$. Number the facets of the container of $V$ so that its
contact labels in facet order are $(0,e_1,e_2,e_3)$. Enumerate $w\in U_1\setminus V$ with
$1\le|w_3|\le5$, all three differences $w-F$ in $U_1$, and
$\conv(F\cup\{w\})$ containing no lattice points other than its vertices. For each
$j\in\{1,2,3\}$, assign $w$ to facet $j$ of the new container and assign
the old contact label on facet $j$ to facet four. The other two labels
retain their facets by Lemma~\ref{lem:transfer}. This gives $3657$ pairs
of adjacent cells with facet assignments. For each pair, try each of
its index-one cells as the initial cell and all common permutations of
the four facets. Translate its first contact label to zero and use the
other three translated labels as a basis. Translate
the other cell so that it agrees on $F$, whenever such a common face
and contact transfer exist. Choose the least resulting $(w,j)$
lexicographically. Each change used here is unimodular and preserves
all incidences and strict comparisons. There are $552$ representatives.
The four individual strict row tests of Section~\ref{sec:checks} eliminate $170$ of these. For each prescribed cell, translate its least vertex to zero and
test~\eqref{eq:or} at each point of $U_1$ outside its vertex set. Whenever only one facet can satisfy the exit condition at a point,
impose that inequality and repeat. These tests eliminate $184$ more pairs.

\subsection{Candidate neighbours}
For each of the $198$ surviving pairs, generate the initial list of
possible cells and facet assignments as follows. The
first cell is the coordinate simplex. In its lattice basis, the set
of candidate edge labels is $U=U_1=11(\Delta-\Delta)\cap\Z^3$.
Intersect the four row cones with the prescribed contact comparisons
and propagate the exit condition~\eqref{eq:or} for each of the two
cells $P$ at every point of $U\setminus P$, normalizing the least
vertex of each cell to zero throughout.
Let $E$ be the six edge differences of each prescribed cell, taken together.
For each primitive lexicographically positive $a\in U$, require
$|\det(u,v,a)|\le11$ for $u,v\in E$, and
$6V(\conv P,\conv Q,[0,a])\le11$ for each pair of prescribed cells $P,Q$, including
repetitions. Form $(0,a,b,c)$ with $a<b<c$ and all six forward differences
on this list. For the candidate vertex set $C=\{0,a,b,c\}$ require that its closed hull
contain no lattice points other than its four vertices and that its index
lie in $[1,5]$. For prescribed cells $P,Q$ also require
$6V(\conv P,\conv Q,\conv C),6V(\conv P,\conv C,\conv C)\le11$. These tests use the same lattice normalization as
Lemma~\ref{lem:mix}. For three lattice polytopes, $6V$ is the
alternating sum of the ordinary volumes of their eight Minkowski
sums, with a positive sign on the sum of all three. Equivalently,
when the sum of the first two polytopes has full dimension, take
each of its primitive outer facet normals, multiply the normalized mixed area of their exposed faces
by the support value of the third polytope and sum the products.
These are Proposition~2.1(5) and formula~(2.10) of~\cite{ABS},
with normalized mixed volume $6V$. Normals of zero mixed area contribute zero. The alternating-sum
formula also applies when some of the polytopes have smaller dimension.

Try every assignment of facets to vertices, retaining the given assignment
on each prescribed vertex set. Intersect the respective row cones with the comparisons
of the trial cell and require a strict solution in every row. For each
$p\in C$ and each integer point $x\in(p+[-1,1]^3)\setminus C$,
some row satisfying the trial comparisons must place $x$
strictly beyond its facet, as in~\eqref{eq:or}. Impose the contact transfers at
all shared prescribed faces. At a face with two prescribed incidences
discard every third one. In each round delete each candidate facet assignment for which some face
has no possible partner in the current list. Repeat
until no further assignment is removed. For $125$ of the $198$ initial
pairs, a prescribed cell loses its assignment. The remaining $73$ pairs
give the starting points for the search below.

\subsection{Completing the face pairings}
At each split carry the row cones and surviving list forward. Add the
new strict contact comparisons. Continue to propagate the exit condition on the full
window $U$ for the initial cells. For every prescribed cell $P$,
including each newly added one, propagate the same condition at
the integer points in the boxes $p+[-2,2]^3$ about its vertices and at
the eight vertices of the unit cube with lower corner $\lfloor(x+y)/2\rfloor$
for each pair $x,y\in P$. Also test $v$ and $v\pm e_j$, $1\le j\le3$,
for every vertex $v$ of the whole prescribed list. Omit $P$'s own vertices.
Apply these tests to the initial pair as well.
For a candidate $C$ and the new prescribed cell $Q$,
test $6V(\conv C,\conv C,\conv Q)$ and $6V(\conv C,\conv Q,\conv P)$ against eleven for all prescribed $P$
(including $Q$). Keep a possible contact $w_i$ only if the $i$th cone admits
$g_i(w_i-v)>0$ for each other vertex $v$ separately. Again delete any assignment that has lost its last possible partner at some face,
repeating until no further deletion occurs. When at most $1400$ vertex sets remain, use the lifted
edge link test below, alternating the link and face tests until neither removes
another candidate. Test whether~\eqref{eq:cup} has a solution with
$x_P=1$ for every prescribed cell $P$. Disappearance of a prescribed cell or
failure of one of these necessary tests closes the node. Otherwise split
on an unpaired face with the fewest allowed opposite incidences, breaking ties by its lexicographically normalized triangle.
Count every attempted child, including an immediate rejection.

For the \reviewterm{link test} fix a lifted edge $[0,v]$, choosing $v$ positive
lexicographically. Each remaining vertex set with an edge translating to
$[0,v]$ contributes an edge to the link. Translate the cell to align
these edges, then join its other two vertex labels. By
Lemma~\ref{lem:topology} the actual link is one simple cycle. Delete
bridges, since deleting a bridge disconnects its endpoints and no cycle
can use it. Every simple cycle lies in one biconnected block, a maximal
connected subgraph with no cut vertex. Thus all prescribed edges must
lie in one such block. Once they determine the block, remove edges outside it.
Reject prescribed degree above two. At a vertex of prescribed degree
two remove all its other edges. If the prescribed edges contain a
cycle, that cycle must contain all prescribed edges of the link and no other
edge can be used. A deletion removes the corresponding vertex set
and all its facet assignments. These tests concern the lifted edge, so vertices in
different translates are distinguished.

Each actual incidence survives the face test, and the lifted link of
an actual edge survives the link test. Thus at every split an actual
covering gives a child containing its adjacent cell. The $73$ trees have $6709$ states. Of these, $6077$ are rejected and
$612$ split. The remaining $20$ have all faces paired. Their signed index
sum is six in sixteen cases and minus six in four.

\subsection{Normalization and determinant bounds}
For each paired list, try every index-one cell and all $24$ common
permutations of the facets. In the chosen cell translate the label on
facet one to zero. Its other three translated labels, in facet order,
give a lattice basis for the whole list. Now translate the first contact
label of each cell to zero and write its four labels in facet order in
this common basis. Sort the
list of these tuples lexicographically. Minimize the largest absolute
coordinate first and then the sorted list lexicographically.
This normalization gives one representative for each of the supports
I, II and III in Section~\ref{sec:solid-supports}.

For each support start its cones anew from the contact comparisons and
impose the forced exit inequalities for each cell separately, using its
radius-two vertex boxes and the eight cube vertices for each midpoint
as in Section~\ref{sec:checks}. For $\mathcal L$ take $1$,
the forms $L_P=1-\sum_i g_i(w_i)$ and the primitive inward facet normals
of the resulting row cones, expressed as linear forms in $B$. Use the
coefficient inequalities of Section~\ref{sec:last} with products of
three factors and $LH^2$. Here $H$ ranges over homogeneous forms in the
nine entries with at most four nonzero coefficients, each $\pm1$. For each normalized support, \texttt{python3 bounds/run.py solid} checks
the coefficient inequalities for the two signs. The table gives upper
bounds for $-\det B$ and $\det B$ verified by the supplied weights.
For the two bounds zero, every coefficient of $e\det B+P(B)$ vanishes.
\[
\begin{array}{c|r|r|r|rr}
 \text{candidate}&\text{cells}&\text{multiplicity}&|\mathcal L|&
 10^{12}(-\det B)\le&10^{12}\det B\le \\\hline
 \mathrm{I}&6&5&19&78605674800&0\\
 \mathrm{II}&6&1&19&83153893839&0\\
 \mathrm{III}&14&14&27&327&83395943758
\end{array}
\]
The multiplicities sum to twenty. The largest of the six displayed
right sides is $83395943758$, so $10^6D\le83396$ on every retained
all-solid support.

\section{Continuation from a longest flat component}\label{sec:flat-runs}
\subsection{Completion from a path and its initial cap}\label{sec:flat-pairings}
Choose a flat component of maximum length $m$, counting its flat cells.
Use the three edges at one vertex of the initial cap as a basis, so that
the cap has vertex set $\{0,e_1,e_2,e_3\}$.
Prescribe this cap and the successive parallelograms of the word as in
Section~\ref{sec:flat-components}. For the cap, assignment 1 lists the
contact labels as $(e_1,e_3,e_2,0)$ in facet order, and assignment 2
lists them as $(e_2,e_1,e_3,0)$. The empty word requires only assignment 1
by symmetry. Thus the ten words requiring completion give nineteen starting lists. The continuation
uses the face and link tests of Appendix~\ref{sec:unit-runs} with
three changes.
\begin{enumerate}
\item Replace its initial window for edges and exits by
\[
 U=11(\Delta-\Delta)\cap\{(x,y,z)\in\Z^3:
                     |x|+|y|\le11,\ |z|\le2\}.
\]
The index-one cap supplies the first bound, and the first
parallelogram supplies the other two.
\item Admit solid candidates of index $1$ through $8$ and flat
fundamental parallelograms. For a sorted tuple $(0,a,b,c)$ the latter
test asks whether one of $a,b,c$ is the sum of the other two and
whether their cross product has relatively prime coordinates.
Include the mixed-volume tests $6V(\conv P,\conv P',\conv C)\le11$ for every candidate
$C$ and prescribed $P,P'$. The repeated-candidate test
$6V(\conv C,\conv C,\conv P)\le11$ is used when $C$ is solid.
\item After every addition join prescribed flat cells sharing a
triangle up to translation. Reject any component with more than $m$
cells. The remaining row, opposite-incidence, lifted-link and degree
tests, and the rule for breaking ties when choosing a face, are as in
Appendix~\ref{sec:unit-runs}. In the link of $[0,v]$, include one edge for
each edge of a cell that translates to $[0,v]$. For example, both
$[0,a]$ and $[b,a+b]$ of the flat cell $\{0,a,b,a+b\}$ translate to
$[0,a]$, giving the two edges of the link $[b,a+b]$ and $[-b,a-b]$.
When all faces are paired, retain the lists whose signed index sum is $\pm6$.
\end{enumerate}

In the recursion below, $S$ is the list of prescribed cells, $A$ records
the surviving facet assignments on each candidate vertex set, and $C$ is
the tuple of four row cones. Each child inherits its own copy of
$A$ and $C$.
\begin{quote}\small\ttfamily\raggedright
Complete($S,A,C$)\\
\quad Impose the comparisons and all forced exit inequalities for $S$ on $C$.\\
\quad Filter $A$ by mixed volumes, possible facet contacts, and the fixed faces.\\
\quad Reject failure of strict feasibility, excess flat length, or a lost prescribed cell.\\
\quad Propagate face and link exclusions to stability. Reject a lost prescribed cell or failed parity test.\\
\quad If every triangle of $S$ is paired, retain $S$ if its degree is $\pm1$, and return.\\
\quad Choose an unpaired face of $S$ with the fewest partners, using the tie rule above.\\
\quad For each listed contact tuple $P$ call Complete($S\cup\{P\},A,C$).
\end{quote}
Each recursive call adds a new vertex set from a fixed finite list.
Closure of a prescribed connected selection in the actual dual
graph implies that it is the whole quotient complex. It follows by
induction that some retained support represents any covering with the
chosen longest word. For each assignment on the initial cap, the table
counts the states visited in the search and the fully paired lists retained, before
choosing representatives up to change of basis and permutation of the facets.
\[
 \begin{array}{c|rr|rr}
 &\multicolumn{2}{c|}{\text{facet assignment 1}}
 &\multicolumn{2}{c}{\text{facet assignment 2}}\\
 \text{word}&\text{states}&\text{paired lists}&\text{states}&\text{paired lists}\\ \hline
 \varnothing&2387&164&&\\
 \mathsf{L}&2994&436&1859&221\\
 \mathsf{LL}&3128&479&2141&201\\ \hline
 \mathsf{LR}&2967&502&2313&393\\
 \mathsf{LLL}&2855&468&1628&203\\
 \mathsf{LLR}&1050&220&1114&220\\
 \mathsf{LRL}&537&93&527&93\\
 \mathsf{LLLL}&1948&313&1389&165\\
 \mathsf{LLLR}&181&16&155&16\\
 \mathsf{LLLLL}&1267&193&1029&126
 \end{array}
\]

The command \texttt{python3 flat/run.py} runs the recursion on the nineteen
starting lists and normalizes the outputs. The upper panel gives $1501$
paired lists and the lower panel $3021$.
The normalization in Appendix~\ref{sec:unit-runs} gives $399$ distinct
supports in the first group and $945$ in the second.
Number the representatives separately in each group, starting at one.
Order them first by the number of cells and then lexicographically by the
sorted list of contact tuples, with each tuple in facet order. These
numbers are the line numbers in \texttt{flat/short.jsonl.gz} for the $399$
supports and in \texttt{flat/long.jsonl.gz} for the $945$, after decompression.

\subsection{Finding and checking the three-row bounds}\label{sec:ray-weights}
For the affine forms involving $m(v)$ one may compute the mixed
volumes entirely in two dimensions. Choose $k$ with $v_k\ne0$ and,
for the two other indices $i,j$, project $x$ to
$(v_kx_i-v_ix_k,v_kx_j-v_jx_k)$. If $A,A'$ are the projections of
$\conv P,\conv P'$, then
\[
 6V(\conv P,\conv P',[0,v])=
 \frac{\operatorname{area}(A+A')-\operatorname{area}(A)
              -\operatorname{area}(A')}{|v_k|}.
\]
Indeed, adding $[0,v]$ to a polytope increases its volume by the
area of its image under this projection divided by $|v_k|$.
Polarization gives the displayed formula, including degenerate images.

Write the chosen affine forms as $a_s+\sum_j g_j(b_{sj})\ge0$ with
$b_{s4}=0$. If row $t$ is omitted and $p_s\ge0$, put
\[
 w=\sum_s p_sa_s,\qquad
 h_j=\sum_s p_s(b_{st}-b_{sj})\quad(j\ne t).
\]
The relation $\sum_i g_i=0$ rewrites their weighted sum as
$w-\sum_{j\ne t}g_j(h_j)\ge0$, the linear inequality in
Lemma~\ref{lem:three-row}. For the bound $b=D_0$, normalize $w=1$
and try each omitted row.
Require the three ray factors to be nonnegative. For a ray triple
of nonzero determinant put $m=|\det(r,s,u)|/(27b)$. Its factors
$x,y,z$ satisfy $xyz\ge m$ precisely when there exist $U,V\ge0$
such that
\[
                  U^2\le xy,\qquad mV^2\le z,\qquad UV\ge1 .
\]
Indeed these inequalities imply $1\le U^2V^2\le xyz/m$, and the
converse follows by taking $U=\sqrt{xy}$ and $V=1/U$.
The constraints are convex quadratic cone constraints with rational
data. Rounding or rational reconstruction gives trial weights.
Clear denominators and check \eqref{eq:three-rays} on every triple
of primitive rays. A homogeneous form nonnegative on a row cone is
a nonnegative combination of its inward facet forms.
The ancillary files \texttt{long-weights.jsonl.gz} and
\texttt{short-weights.jsonl.gz} in \texttt{bounds/} supply the integer
weights used below. The program \texttt{bounds/run.py} reconstructs the
row cones from the contact tuples and checks each ray inequality.

\subsection{Three-row bounds for the 945 supports}
The seven words in the lower panel of the table in Appendix~\ref{sec:flat-pairings}
give $945$ distinct supports after the normalization of Appendix~\ref{sec:unit-runs}.
For each take the constant $1$, the forms $L_P$ and the primitive
inward facet forms of the tightened row cones. For every representative
except $65$, the supplied weights verify~\eqref{eq:three-rays} directly.
For representative $65$, split the domain along $g_3(0,2,3)=0$ in the
normalized coordinates and adjoin the corresponding signed cut form
on each side. Choosing the omitted row separately for each of the $946$
domains gives $27276$ ray triples in total. In each domain
at most ten positive integer weights, each at most $10^6$, suffice
in~\eqref{eq:three-rays} with $b=D_0$. These bounds allow a finite
search in the original scaling of the forms. Alternatively one can
find weights with $w=1$ as above, clear denominators and check the
ray inequalities with integers. The command \texttt{python3 bounds/run.py long}
checks the supplied weights on the normalized supports.

\subsection{The empty word, L and LL}
The upper panel of the table in Appendix~\ref{sec:flat-pairings} gives five
pairs of a word and an assignment on its initial cap. Together they give
$399$ distinct supports. A support is counted once in each row in which
it occurs, so the rows can overlap. The last four columns count supports
according to the test used. The linear column records a contradiction,
the next two use the three-row bound~\eqref{eq:three-rays} directly or after
cuts, and the last uses a cubic bound.
For the negative sign at representative $104$ this last step uses the
identity in Appendix~\ref{sec:sharp-identity}. The command
\texttt{python3 bounds/run.py short} checks the bounds proved in this subsection.
\[
\begin{array}{cc|r|rrrr}
 & & & & \multicolumn{2}{c}{\text{three-row bound}} &\\
 \text{word}&\text{assignment}&\text{supports}&\text{linear}
 &\text{direct}&\text{with cuts}&\text{cubic}\\\hline
 \varnothing&1&39&0&18&12&9\\
 \mathsf{L}&1&173&0&149&24&0\\
 \mathsf{L}&2&119&0&99&20&0\\
 \mathsf{LL}&1&184&2&177&4&1\\
 \mathsf{LL}&2&115&0&111&4&0
\end{array}
\]
For representatives $341$ and $388$ write
$L_j$ for $L_P$ at the $j$th cell in the sorted list, and $M_{x,y,z}$
for $12\gamma(Bv)-m(v)$ at $v=(x,y,z)$. For $341$ we obtain
\[
\begin{aligned}
 -4={}&720L_7+660L_{12}+12L_{13}+360L_{14}+192g_1(0,-2,1)\\
     &+13M_{1,2,-1}+68M_{2,1,-2}+18M_{2,2,-1}
       +54M_{3,-1,0}+19M_{3,1,-2}.
\end{aligned}
\]
Here $g_1(0,-2,1)$ is an inward facet form of the first row cone.
For $388$ the identity is
\[
\begin{aligned}
 -24={}&204L_1+132L_5+156L_8+228L_{11}+144L_{14}\\
       &+17M_{1,-1,-1}+19M_{1,-1,0}+23M_{1,1,-1}+25M_{1,1,4}.
\end{aligned}
\]
The right sides are nonnegative, excluding both representatives.
For $347$ representatives the three-row bound~\eqref{eq:three-rays}
gives $D\le D_0$ directly.

For a further $40$ representatives the same test is used on subdomains. Consider
primitive $0<n$ lexicographically with $\|n\|_\infty\le3$. In each
row cone put $z_r=(r\cdot n)/\|r\|_1$ on its primitive rays. For
$\min z_r<0<\max z_r$ rank the possible cuts by decreasing
$(-\min z_r)(\max z_r)/\|n\|_1^2$, breaking ties by decreasing row
number and then decreasing lexicographic $n$. Try the first three cuts
at depths zero and one and the first cut at depths two, three and four.
Intersect that row cone with the two halfspaces and adjoin $g_i(n)$
or $-g_i(n)$ to the respective affine lists. A cut succeeds only when
both children succeed. The successful trees have $191$ terminal
subdomains in total and depth at most five. On these and the $347$
unsplit domains, \eqref{eq:three-rays} holds with $b=D_0$ using at most
$30$ positive integer weights per domain, each below $10^9$.

The ten remaining representatives are
$1,2,3,4,8,27,42,97,104,309$. For these supports $\mathcal L_0$ consists of the constant $1$, the forms $L_P$
and the primitive facet forms of the four tightened cones. For each
determinant sign $e$ use all
triple products drawn with repetition from $\mathcal L_0$. Also allow
$LH^2$, $L\in\mathcal L_0$, with $H$ a homogeneous form in the nine
entries having at most four nonzero coefficients, all $\pm1$. There
are $1425$ possible $H$ up to sign. Solve the rational linear
feasibility problem in~\eqref{eq:cubic-terminal} with right side
$83562/10^6$ for both signs on the nine representatives other than
$104$. All eighteen systems are feasible. For $104$ the positive sign
also gives a feasible system with right side $17708/10^6$. For the
negative sign we need the sharp bound and its equality case.

\newpage
\subsection{The sharp determinant identity}\label{sec:sharp-identity}
The thirteen contact tuples of short-word representative $104$ specify the
cells of the $84$-shape cycle in normalized lattice coordinates. Their
basis matrix $B$ can still vary. For the remaining sign, $\det B<0$, the
covolume is $D=-\det B$ and the density is $1/(6D)$. The $84$-shape has
covolume $21/250$. To prove that its density $125/63$ is optimal on this
support, we express the difference in covolumes as
\[
       \frac{21}{250}-D=\frac{21}{250}+\det B
         =\sum_{j=1}^{215}\lambda_j Q_j(B),\qquad \lambda_j>0 .
\]
The containment inequalities and the row cones of this support make each
product $Q_j(B)$ nonnegative on the domain and, by continuity, on its closure.
The identity thus gives
$D\le21/250$. At equality all the products must vanish. We use their
common zeros to determine the unique normalized matrix $B_*$.

\emph{The product identity.}
Write $a=e_1,b=e_2,c=e_3$ and number the nonzero labels by
\[
\begin{array}{r|rrrrrrr}
j&1&2&3&4&5&6&7\\
 v_j&c-3a&c-2a-b&c-2a&c-a-b&c-a&b-a&c-b\\ \hline
 j&8&9&10&11&12&13&14\\
 v_j&c&b&a-b+c&a&a+b&2a-b&3a+b-c
\end{array}
\]
In every cell the label on facet one is zero. The triples below give
the subscripts $j$ of the labels $v_j$ on facets two, three and four,
respectively, in each of the thirteen cells.
\[
\begin{gathered}
 (2,3,1),\ (2,5,3),\ (4,11,2),\ (7,11,4),\ (8,9,6),\ (10,11,7),\\
 (11,5,2),\ (11,6,5),\ (11,8,6),\ (11,9,8),\
 (11,14,9),\ (13,10,7),\ (14,12,9).
\end{gathered}
\]
Prepending zero to each triple gives the vertex order of the corresponding
oriented simplex in~\eqref{eq:84-cycle}, here written in integer labels. The points $u_j$ there
are $Av_j$, where $A$ has columns $-C_1,C_3-C_1,C_3-C_2-3C_1$
and $C_i$ are the columns of~\eqref{eq:84-lattice}.

Each cell $P$ supplies a nonnegative factor $L_P(B)=1-r_B(P)$ by
containment. Use these thirteen factors and the constant $1$.
The other sixteen factors are the inward facet forms $g_i(n)$ of the row cones,
as $n$ ranges over the four normals in row $i$ of the following table.
\[
\begin{array}{c|rrrr}
 i&\multicolumn{4}{c}{n}\\\hline
 1&(-3,-1,1)&(-2,1,0)&(1,-1,0)&(2,1,-1)\\
 2&(-2,-1,1)&(0,0,1)&(1,0,-1)&(3,1,-1)\\
 3&(-2,0,1)&(0,1,-1)&(1,-1,1)&(3,0,-1)\\
 4&(0,1,0)&(0,-1,1)&(-1,1,-1)&(-1,-1,0)
\end{array}
\]
As in~\eqref{eq:cubic-terminal}, we use products of three of these
factors and products $LH^2$, with $L$ one of the factors and $H$ one of
the homogeneous forms specified in the preceding subsection. Each product
is nonnegative on the domain. To identify the $215$ products used in the
identity, sort the thirty factors lexicographically by their
ten coefficients, constant first and entries of $B$ in row order.
Form the coefficient columns of their $\binom{32}{3}$ triple products
in lexicographic order of nondecreasing factor triples. Append the
columns for $LH^2$ in lexicographic order of $(L,H)$, with the first
nonzero coefficient of $H$ positive. Number columns starting at one.

\begin{samepage}
The integers below give the successive differences between the selected
column numbers, starting from zero. Taking cumulative sums recovers the
columns of $Q_1,\ldots,Q_{215}$. For example, the first four entries
$46,2,31,3$ select columns $46,48,79,82$.
{\small
\begin{verbatim}
46 2 31 3 19 26 78 3 26 13 3 21 15 3 10 8 11 9 28 5 11 29 26 198 23
21 56 4 19 7 14 149 9 20 3 6 15 7 23 2 2 26 45 11 1 48 12 6 29 9
11 5 30 25 38 20 2 3 21 2 3 22 22 8 5 16 28 11 6 19 2 16 2 5 18
36 5 7 36 12 72 60 16 1 2 42 15 4 2 10 18 5 4 9 23 8 6 1 41 12
95 11 17 16 1 11 11 5 3 10 2 4 58 6 15 2 8 26 12 86 1 1 43 1 1
52 7 1 1 16 11 1 14 8 16 4 48 52 122 12 9 48 77 31 13 1 10 34 18 11
1 14 22 127 15 20 15 3 16 8 13 8 18 30 116 4 2 11 20 11 8 8 75 16 14
5 22 9 20 1 2 16 19 77 7 20 9 1 9 14 13 2 44 8 44 28 8 13 13 5
55 55 14 21 25 47 21 10 3 80 465 5 4 37459 517
\end{verbatim}
}
\end{samepage}
To find the weights, take these $215$ columns as the matrix $M$,
with one row for each monomial of degree at most three in the nine
entries of $B$. Let $b$ be the coefficient vector of $21/250+\det B$
in the same order. Matching coefficients gives $M\lambda=b$.
Rational elimination gives a solution with every $\lambda_j>0$.
The rank of $M$ modulo $1009$ is $215$, so the weights are unique.
For example, one can select independent rows by elimination modulo
$1009$, solve those $215$ equations over $\mathbb Q$, and substitute
in all $220$ coefficient equations and check that every entry is
positive. This proves the identity. The normalized matrix of the
$84$-shape is
\[
           B_*=\tfrac{1}{10}\begin{pmatrix}-2&-3&-8\\
                                    2&-4&1\\2&6&5\end{pmatrix}.
\]
At $B_*$ the thirteen forms $L_P$ vanish, the sixteen homogeneous factors
equal $1/10$, and the determinant is $-21/250$.

\emph{Equality.}
Suppose $-\det B=21/250$ in the closure of the domain of this support.
All $215$ weights are positive, so every product in the identity
vanishes. One of them eliminates a matrix entry. The last selected
column, $42826$, is the square of
\[
                H=B_{12}+B_{13}+B_{32}+B_{33}.
\]
Thus $H=0$. It remains to see why the eight other entries must agree
with those of $B_*$. For this it suffices to use the $213$ selected
triple products. Substitute $B_{33}=-B_{12}-B_{13}-B_{32}$ in them.
List the eight remaining entries in row
order as $x_1,\ldots,x_8$, put $x_0=1$, and expand each product in
the monomials $x_i x_j x_k$, $0\le i\le j\le k\le8$. Row reduction
of the resulting $213\times165$ integer matrix modulo $1009$ gives
$164$ pivots.

Over $\mathbb Q$ these products therefore span a space of dimension
at least $164$. They all vanish at $B_*$. Evaluation there has
a $164$-dimensional kernel in the $165$-dimensional polynomial
space, so the products span this kernel. In particular each of the
eight linear coordinate differences $x_i-x_i(B_*)$ belongs to their
span. Thus equality forces the eight entries to agree with $B_*$,
and $H=0$ gives the ninth. This proves uniqueness on the closed domain. The command
\texttt{python3 bounds/run.py sharp} verifies the coefficient identity and the
ranks used for the equality case.

\end{document}